\documentclass{article}
\usepackage{cite}
\usepackage{graphicx,epsfig}
\usepackage{amsmath,amsthm}
\usepackage{subfigure}
\usepackage{array}
\usepackage{indentfirst}
\usepackage{amssymb}
\usepackage{setspace}
\fontsize{12}{12pt}
\theoremstyle{plain}
\newtheorem{theo}{Theorem}[section]

\newtheorem{prop}[theo]{Proposition}
\theoremstyle{definition}
\newtheorem{defi}{Definition}[section]
\theoremstyle{remark}
\newtheorem{rem}{Remark}[section]
\numberwithin{equation}{section}

\newcommand{\coor}{q_1,\cdots,q_n,p_1,\cdots,p_n}
\newcommand{\Coor}{Q_1,\cdots,Q_n,P_1,\cdots,P_n}
\newcommand{\coorfjpr}{(q_{I_p},p_{\bar{I}_p},Q_{K_r},P_{\bar{K}_r},t)}
\newcommand{\coorfjpri}{(q_{I_p},p_{\bar{I}_p},q_{0_{K_r}},p_{0_{\bar{K}_r}},t),t)}

\newcommand{\pfrac}[2]{\frac{\partial #1}{\partial #2}}
\newcommand{\ppfrac}[3]{\frac{\partial^2 #1}{\partial #2 \partial #3}}

\newcommand{\undemi}{\frac{1}{2}}
\newenvironment{disarray}
 {\everymath{\displaystyle\everymath{}}\array}
 {\endarray}
\newcommand{\hq}{ \Delta q}
\newcommand{\hp}{\Delta p}

\newcommand{\kqz}{\Delta q_0}
\newcommand{\kpz}{\Delta p_0}

\newcommand{\titleone}{\emph{Search in time domain}}

\newcommand{\titlethree}{\emph{Search in position space}}

\begin{document}
\title{Solving two-point boundary value problems using generating functions: Theory and Applications to optimal control
and the study of Hamiltonian dynamical systems\footnote{Abbreviated title: Solving two-point
boundary value problems.}}
\author{
\begin{tabular}{c}
V.M. Guibout\footnote{Graduate Research Assistant, PhD Candidate, Aerospace Engineering Department,
FXB Building, 1320 Beal Avenue, Ann Arbor, MI 48109-2140, guibout@umich.edu } \hspace{2pt} and
D.J.Scheeres
\footnote{Associate Professor, Senior member AIAA, Aerospace Engineering Department, FXB Building,
1320 Beal Avenue, Ann Arbor, MI 48109-2140, scheeres@umich.edu}\\
\textit{University of Michigan, Ann Arbor, Michigan}
\end{tabular}}
\date{}
\maketitle
\begin{abstract}
A methodology for solving two-point boundary value problems in phase space for Hamiltonian systems
is presented. Using Hamilton-Jacobi theory in conjunction with the canonical transformation induced
by the phase flow, we show that the generating functions for this transformation solve any
two-point boundary value problem in phase space. Properties of the generating functions are
exposed, we especially emphasize multiple solutions, singularities, relations with the state
transition matrix and symmetries. Then, we show that using Hamilton's principal function we are
also able to solve two-point boundary value problems, nevertheless both methodologies have
fundamental differences that we explore. Finally, we present some applications of this theory.
Using the generating functions for the phase flow canonical transformation we are able to solve the
optimal control problem (without an initial guess), to study phase space structures in Hamiltonian
dynamical systems (periodic orbits, equilibrium points) and classical targeting problems (this last
topic finds its applications in the design of spacecraft formation trajectories, reconfiguration,
formation keeping, etc...).

\end{abstract}

\section{Introduction}
One of the most famous two-point boundary value problems in astrodynamics is Lambert's problem,
which consists of finding a trajectory in the two-body problem which goes through two given points
in a given lapse of time. Even though the two-body problem is integrable, no analytical solution
has been found to this problem so far, and solving Lambert's problem still requires one to solve
Kepler's equation, which has motivated many papers since $1650$\cite{col93}. For a general
Hamiltonian dynamical system, a two-point boundary value problem is solved using shooting methods
combined with Newton iteration. Though very systematic, this technique requires a ``good" initial
guess for convergence and is not appropriate when several boundary value problems need to be
solved. In order to design a change of configuration of a formation of $n$ spacecraft, $n!$
two-point boundary value problems need to be solved\cite{wan99had}, hence for a large collection of
spacecraft the shooting method is not efficient. In this paper we address a technique which allows
us to solve $m$ boundary value problems at the cost of $m$ function evaluations once generating
functions for the canonical transformation induced by the phase flow are known. These generating
functions are solutions of the Hamilton-Jacobi equation and for a certain class of problem they can
be found offline, that is during mission planning. Moreover, the theory we expose allows us to
formally solve any kind of two-point boundary value problem, that is, given a $n$-dimensional
Hamiltonian system and $2n$ coordinates among the $4n$ defining two points in the phase space, we
find the other $2n$ coordinates. The Lambert problem is a particular case of this problem where the
dynamics is Keplerian, the position of two points are given and the corresponding momenta need to
be found. Another instance of such a problem is the search for trajectories which go through two
given points in the momentum space (i.e., the conjugate of the Lambert problem). Properties of the
solutions found are studied, such as multiple solutions, symmetries and relation to the state
transition matrix for linear systems. Then, we expose another method to solve two-point boundary
value problems based on Hamilton's principal function and study how it compares to generating
functions. Finally, we present direct applications of this theory through the optimal control
problem and the study of some Hamiltonian dynamical systems. Solving the optimal control problem
using generating functions was first introduced by Scheeres et al.\cite{sch03gui}, we will review
their method in this paper and expand it to more general optimal control problems. Applications to
Hamiltonian dynamical systems were first studied by Guibout and Scheeres\cite{gui02,gui03b} for
spacecraft formation flight design and for the computation of periodic orbits.

\section{Solving a two-point boundary value problem}
In this section, we recall the principle of least action for Hamiltonian systems and derive the
Hamilton-Jacobi equation. Local existence of generating functions is proved. We underline that we
do not study global properties. In general, we do not know a priori if the generating functions
will be defined for all time and in most of the cases we found that they develop singularities. We
refer the reader to \cite{abr78,arn88a,gol65,gre77,gui02,lan77,mar98} for more details on local
Hamilton-Jacobi theory, \cite{abr78,arn88a,mar98} for global theory and
\cite{ehl00new,abr78,arn88a} and section \ref{sect:singularities} of this paper for a study of
singularities.

\subsection{The Hamilton-Jacobi theory}

Let $(P,\omega,X_H)$ be a Hamiltonian system with $n$ degrees of freedom, and $H:P \times
\mathbb{R}\rightarrow
\mathbb{R}$ the Hamiltonian function. In the extended phase space $P \times \mathbb{R}$, we consider an integral curve
of the vector field $X_H$ connecting the points $(q_0,p_0,t_0)$ and $(q_1,p_1,t_1)$. The principle
of least action reads:

\begin{theo}(\textbf{The principle of least action in phase space})
The integral $\int_0^1 pdq-Hdt$ has an extremal in the class of curve $\gamma$ whose ends lie in
the n-dimensional subspaces $(t=t_0,q=q_0)$ and $(t=t_1,q=q_1)$ of extended phase space.
\end{theo}
\begin{proof}
We proceed to the computation of the variation.
\begin{eqnarray}
\delta\int_\gamma (p\dot{q}-H)dt&=&\int_\gamma \left(\dot{q}\delta p+p \delta
\dot{q}-\pfrac{H}{q}\delta q- \pfrac{H}{p}\delta p \right)dt\nonumber\\
&=&\left[ p\delta q \right]_0^1+\int_\gamma \left[ \left(\dot{q}-\pfrac{H}{p}\right)\delta p
-\left(\dot{p}+\pfrac{H}{q}\right)\delta q \right]dt
\end{eqnarray}
Therefore, since the variation vanishes at the end points, the integral curves of the Hamiltonian
vector field are the only extremals.
\end{proof}

\begin{rem}
The condition for a curve $\gamma$ to be an extremal of a functional does not depend on the choice
of coordinate system, therefore the principle of least action is coordinate invariant.
\end{rem}

Now let $(P_1, \omega_1)$ and $(P_2,\omega_2)$ be symplectic manifolds,
\begin{defi}
A smooth map $f:P_1\times \mathbb{R} \rightarrow  P_2\times \mathbb{R}$ is a canonical
transformation if and only if

(1)- $f$ is a $C^\infty$-diffeomorphism,

(2)- $f$ preserves the time, i.e., there exists a function $g_t$ such that $f(x,t)=(g_t(x),t)$,

(3)- for each $t$, $g_t: P_1\rightarrow P_2$ as defined above is a symplectic diffeomorphism and
$f$ preserves the canonical form of Hamilton's equations.
\end{defi}
All three points in this definition are not independent but we mention them for sake of clarity. It
can be proved\cite{abr78} that if $g_t$ is symplectic then $f$ is a diffeomorphism. Moreover, the
third point of the definition differs from book to book. We chose Abraham's definition\cite{abr78}
but very often the third item reduces to ``$f$ preserves Hamilton's equations" (Goldstein
\cite{gol65}, Greenwood \cite{gre77}). Arnold \cite{arn88a} argues that this definition differs
from the original definition, the third item should actually be ``$g_t$ is symplectic" which
implies, but is not equivalent to, ``$f$ preserves the canonical form of Hamilton's equations".

Consider a canonical transformation $f:(q_i,p_i,t)\mapsto (Q_i,P_i,t)$. Since Hamilton's equations
are preserved, we have:
\begin{equation}
\left\{
\begin{disarray}{lcl}
\dot{Q}_i&=&\pfrac{K}{P_i}\\
&&\\
\dot{P}_i&=&-\pfrac{K}{Q_i}
\end{disarray}
\right.
\end{equation}
where $K=K(Q,P,t)$ is the Hamiltonian of the system in the new set of coordinates.

On the other hand, we have seen that the principle of least action is coordinate invariant. Hence:
\begin{equation}
\delta \int_{t_0}^{t_1}\left(\sum_{i=1}^np_i
\dot{q}_i-H(q,p,t)\right)dt=0
\label{eq:pla_old}
\end{equation}
\begin{equation}
\delta \int_{t_0}^{t_1}\left(\sum_{i=1}^nP_i
\dot{Q}_i-K(Q,P,t)\right)dt=0
\label{eq:pla_new}
\end{equation}

From Eqns. \ref{eq:pla_old} - \ref{eq:pla_new}, we conclude that the integrands of the two
integrals differ at most by a total time derivative of an arbitrary function $F$:
\begin{equation}
\sum_{i=1}^n p_idq_i-Hdt=\sum_{j=1}^n P_j dQ_j -Kdt + dF
\label{eq:canonical_transf_relation_f1}
\end{equation}

Such a function is called a generating function for the canonical transformation $f$ and is,
\textit{a priori}, a function of both the old and the new variables and time. The two sets of
coordinates being connected by the $2n$ equations, namely, $f(\coor,t)=(\Coor,t)$, $F$ can be
reduced to a function of $2n+1$ variables among the $4 n+1$. Hence, we can define $4^n$ generating
functions that have $n$ variables in $P_1$ and $n$ in $P_2$. Among these are the four kinds defined
by Goldstein\cite{gol65}, $F_1(q_1,\cdots,q_n,Q_1,\cdots,Q_n,t)$,
$F_2(q_1,\cdots,q_n,P_1,\cdots,P_n,t)$, $F_3(p_1,\cdots,p_n,Q_1,\cdots,Q_n,t)$ and \newline
$F_4(p_1,\cdots,p_n,P_1,\cdots,P_n,t)$.

Let us first consider the generating function $F_1(q,Q,t)$. The total time derivative of $F_1$
reads:
\begin{equation}
dF_1(q,Q,t)=\sum_{i=1}^n \pfrac{F_1}{q_i}dq_i+\sum_{j=1}^n
\pfrac{F_1}{Q_i}dQ_i+\pfrac{F_1}{t}dt
\end{equation}
Hence Eq. \ref{eq:canonical_transf_relation_f1} yields:
\begin{equation}
\sum_{i=1}^n (p_i-\pfrac{F_1}{q_i})dq_i-H dt=\sum_{j=1}^n
(P_j+\pfrac{F_1}{Q_j}) dQ_j -K dt+\pfrac{F_1}{t}dt
\label{eq:hamiltonjacobi_f1}
\end{equation}
Assume that $(q,Q,t)$ is a set of independent variables, then Eq. \ref{eq:hamiltonjacobi_f1} is
equivalent to:
\begin{eqnarray}
p_i&=&\pfrac{F_1}{q_i}(q,Q,t)\label{eq:f1_in_generala} \\
P_i&=&-\pfrac{F_1}{Q_i}(q,Q,t) \label{eq:f1_in_generalb}\\
K(Q,-\pfrac{F_1}{Q},t)&=&H(q,\pfrac{F_1}{q},t)+\pfrac{F_1}{t}
\label{eq:f1_in_generalc}
\end{eqnarray}
If $(q,Q)$ is not a set of independent variables, we say that $F_1$ is singular.

Now let us consider more general generating functions. Let $(i_1,\cdots,i_p)(i_{p+1},\cdots,i_{n})$
and
\newline$(k_1,\cdots,k_r)(k_{r+1},\cdots,k_{n})$ be two partitions of the set $(1,\cdots,n)$ into two
non-intersecting parts such that $i_1<\cdots<i_p$, $i_{p+1}<\cdots<i_{n}$, $k_1<\cdots<k_r$ and
$k_{r+1}<\cdots<k_{n}$ and define $I_p=(i_1,\cdots,i_p)$, $\bar{I}_p=(i_{p+1},\cdots,i_{n})$,
$K_r=(k_1,\cdots,k_r)$ and $\bar{K}_r=(k_{r+1},\cdots,k_{n})$. If
$$(q_{I_p},p_{\bar{I}_p},Q_{K_r},P_{\bar{K}_r})=
(q_{i_1},\cdots,q_{i_p},p_{i_{p+1}},\cdots,p_{i_{n}},Q_{k_1},\cdots,Q_{k_r},P_{k_{r+1}},\cdots,P_{k_{n}})$$
are independent variables, then we can define the generating function $F_{I_p,K_r}$:
\begin{equation}
F_{I_p,K_r}(q_{I_p},p_{\bar{I}_p},Q_{K_r},P_{\bar{K}_r},t)=
F(q_{i_1},\cdots,q_{i_p},p_{i_{p+1}},\cdots,p_{i_{n}},Q_{k_1},\cdots,Q_{k_r},P_{k_{r+1}},
\cdots,P_{k_{n}},t)
\end{equation}
Expanding $dF_{I_p,K_r}$ yields:
\begin{equation}
dF_{I_p,K_r}=
\sum_{a=1}^p \pfrac{F_{I_p,K_r}}{q_{i_a}}dq_{i_a}+
\sum_{a=p+1}^{n} \pfrac{F_{I_p,K_r}}{p_{i_a}}dp_{i_a}+
\sum_{a=1}^r \pfrac{F_{I_p,K_r}}{Q_{k_a}}dQ_{k_a}+
\sum_{a=r+1}^{n} \pfrac{F_{I_p,K_r}}{P_{k_a}}dP_{k_a}+
\pfrac{F_{I_p,K_r}}{t}dt
\label{eq:dfjrp}
\end{equation}

 and rewriting Eq. \ref{eq:canonical_transf_relation_f1} as a function of the linearly
 independent variables leads to:
\begin{equation}
\sum_{a=1}^p      p_{i_a} dq_{i_a}-
\sum_{a=p+1}^{n} q_{i_a} dp_{i_a}- Hdt=
\sum_{a=1}^r      P_{k_a} dQ_{k_a}-
\sum_{a=r+1}^{n} Q_{k_a} dP_{k_a}- Kdt
+ dF_{I_p,K_r}
\label{eq:varfjrp}
\end{equation}
where $F_{I_p,K_r}=F_1+\sum_{a=r+1}^{n} Q_{k_a}P_{k_a}-\sum_{a=p+1}^{n} q_{i_a}p_{i_a}$ This last
relation defines the Legendre transformation, which allows one to transform one generating function
into another.

Then Eq. \ref{eq:varfjrp} reads:
\begin{eqnarray}
\sum_{a=1}^r    ( P_{k_a}+ \pfrac{F_{I_p,K_r}}{Q_{k_a}})dQ_{k_a}&-&
\sum_{a=r+1}^{n}(-Q_{k_a}+ \pfrac{F_{I_p,K_r}}{P_{k_a}})dP_{k_a}- K dt
+ \pfrac{F_{I_p,K_r}}{t}dt\nonumber\\&=&
\sum_{a=1}^p    ( p_{i_a}- \pfrac{F_{I_p,K_r}}{q_{i_a}})dq_{i_a}
\sum_{a=p+1}^{n}(-q_{i_a}- \pfrac{F_{I_p,K_r}}{p_{i_a}})dp_{i_a}-H dt\nonumber\\&&
\label{eq:hamiltonjacobi_fjrp}
\end{eqnarray}
which is equivalent to:
\begin{eqnarray}
p_{I_p}&=&\pfrac{F_{I_p,K_r}}{q_{I_p}}\coorfjpr \label{eq:fjpr_in_generala}\\
q_{\bar{I}_p}&=&-\pfrac{F_{I_p,K_r}}{q_{\bar{I}_p}}\coorfjpr\label{eq:fjpr_in_generalb} \\
P_{K_r}&=&-\pfrac{F_{I_p,K_r}}{Q_{K_r}}\coorfjpr\label{eq:fjpr_in_generalc} \\
Q_{\bar{K}_r}&=&\pfrac{F_{I_p,K_r}}{P_{\bar{K}_r}}\coorfjpr \label{eq:fjpr_in_generald}\\
K(Q_{K_r},\pfrac{F_{I_p,K_r}}{P_{\bar{K}_r}},-\pfrac{F_{I_p,K_r}}{Q_{K_r}},P_{\bar{K}_r},t)&=&
H(q_{I_p},-\pfrac{F_{I_p,K_r}}{p_{\bar{I}_p}},\pfrac{F_{I_p,K_r}}{q_{I_p}},p_{\bar{I}_p},t)+\pfrac{F_{I_p,K_r}}{t}
\label{eq:fjpr_in_generale}
\end{eqnarray}

For the case where the partitions are $(1,\cdots,n)()$ and $()(1,\cdots,n)$ (i.e.,$p=n$ and $r=0$),
we recover the generating function $F_2$, which verifies the following equations:
\begin{eqnarray}
p_i&=&\pfrac{F_2}{q_i}(q,P,t)\label{eq:f2_in_generala} \\
Q_i&=&\pfrac{F_2}{P_i}(q,P,t)\label{eq:f2_in_generalb} \\
K(\pfrac{F_2}{P},P,t)&=&H(q,\pfrac{F_2}{q},t)+\pfrac{F_2}{t}
\label{eq:f2_in_generalc}
\end{eqnarray}
The case $p=0$ and $r=n$ corresponds to a generating function of the third kind, $F_3$:
\begin{eqnarray}
q_i&=&-\pfrac{F_3}{p_i}(p,Q,t) \label{eq:f3_in_generala}\\
P_i&=&-\pfrac{F_3}{Q_i}(p,Q,t) \label{eq:f3_in_generalb}\\
K(Q,-\pfrac{F_3}{Q},t)&=&H(-\pfrac{F_3}{p},p,t)+\pfrac{F_3}{t}
\label{eq:f3_in_generalc}
\end{eqnarray}
Finally, if $p=0$ and $r=0$, we obtain $F_4$:
\begin{eqnarray}
q_i&=&\pfrac{F_4}{p_i}(p,P,t) \label{eq:f4_in_generala}\\
Q_i&=&-\pfrac{F_4}{P_i}(p,P,t)\label{eq:f4_in_generalb} \\
K(-\pfrac{F_4}{P},P,t)&=&H(\pfrac{F_4}{p},p,t)+\pfrac{F_4}{t}
\label{eq:f4_in_generalc}
\end{eqnarray}

\subsection{The phase flow is a canonical transformation}
In the following we focus on a specific canonical transformation, the one induced by the phase
flow. Let $\Phi_t$ be the flow of an Hamiltonian system:
\begin{eqnarray}
 \label{eq:def_Phi}\Phi_t
:P&\rightarrow&P\nonumber\\
(q_0,p_0)&\mapsto& (\Phi^1_t(q_0,p_0)=q(q_0,p_0,t),\Phi^2_t(q_0,p_0)=p(q_0,p_0,t))\end{eqnarray}
Then, the phase flow induces a transformation $\phi$ on $P\times \mathbb{R}$ defined as follows:
\begin{equation}
\phi:(q_0,p_0,t)\mapsto (\Phi_t(q_0,p_0),t)
\end{equation}

\begin{theo}
The transformation $\phi$ induced by the phase flow is canonical.
\end{theo}
\begin{proof}
The proof of this theorem can be found in Arnold\cite{arn88a}, it is based on the integral
invariant of Poincar{\'e}-Cartan.
\end{proof}

For such a transformation, $(Q,P)$ represents the initial conditions of the system $(q_0,p_0)$, the
Hamiltonian function $K$ is a constant that can be chosen to be $0$ and the equations verified by
the generating function $F_{I_p,K_r}$ become:
\begin{eqnarray}
p_{I_p}&=&\pfrac{F_{I_p,K_r}}{q_{I_p}}\coorfjpri \\
q_{\bar{I}_p}&=&-\pfrac{F_{I_p,K_r}}{p_{\bar{I}_p}}\coorfjpri \\
p_{0_{K_r}}&=&-\pfrac{F_{I_p,K_r}}{q_{0_{K_r}}}\coorfjpri\\
q_{0_{\bar{K}_r}}&=&\pfrac{F_{I_p,K_r}}{p_{0_{\bar{K}_r}}}\coorfjpri \\
0&=&H(q_{I_p},-\pfrac{F_{I_p,K_r}}{p_{\bar{I}_p}},\pfrac{F_{I_p,K_r}}{q_{I_p}},p_{\bar{I}_p},t)+\pfrac{F_{I_p,K_r}}{t}
\label{eq:fjpr_in_general_hj}
\end{eqnarray}
The last equation is often referred to as the Hamilton-Jacobi equation. To solve this equation, one
needs boundary conditions. At the initial time, position and momentum $(q,p)$ are equal to the
initial conditions $(q_0,p_0)$. Hence, $F_{I_p,K_r}$ must generate the identity transformation at
the initial time.

\subsection{Properties of the canonical transformation induced by the phase flow}
In this section we study the properties of generating functions for the phase flow canonical
transformation. First we show that they solve a two-point boundary value problem, and then we prove
a few results on singularities, symmetries and differentiability. In particular, we relate the
generating functions and the state transition matrix for a linear system.

\subsubsection{Solving a two-point boundary value problem}
Consider two points in phase space, $X_0=(q_0,p_0)$ and $X_1=(q,p)$, and two partitions of
$(1,\cdots,n)$ into two non-intersecting parts, $(i_1,\cdots,i_p)(i_{p+1},\cdots,i_{n})$ and
$(k_1,\cdots,k_r)(k_{r+1},\cdots,k_{n})$. A two-point boundary value problem is formulated as
follows: \newline Given $2n$ coordinates $(q_{i_1}, \cdots, q_{i_p}, p_{i_{p+1}},
\cdots,p_{i_{n}})$ and $(q_{0_{k_1}},\cdots,q_{0_{k_r}},p_{0_{k_{r+1}}},\cdots,p_{0_{k_{n}}})$,
find the remaining $2n$ variables such that a particle starting at $X_0$ will reach $X_1$ in $T$
units of time.

From the relationship defined by Eqns. \ref{eq:fjpr_in_generala}, \ref{eq:fjpr_in_generalb},
\ref{eq:fjpr_in_generalc} and \ref{eq:fjpr_in_generald}, we see that the generating function
$F_{I_p,K_r}$ solves this problem. Lambert's problem is a particular case of boundary value problem
where the partitions of $(1,\cdots,n)$ are $(1,\cdots,n)()$ and $(1,\cdots,n)()$. Though, given two
positions $q_f$ and $q_0$ and a transfer time $T$, the corresponding momentum vectors are found
from the relationships verified by $F_1$:
\begin{eqnarray}
p_i&=&\pfrac{F_1}{q_i}(q,q_0,T) \nonumber\\
p_{0_i}&=&-\pfrac{F_1}{q_{0_i}}(q,q_0,T)
\end{eqnarray}

\subsubsection{Existence and properties of the generating functions}
In the first section we proved the existence of a generating function using the assumption that its
variables are linearly independent. This is not always true at every instant. As an example let us
look at the harmonic oscillator. The equations of motion are given by:
\begin{eqnarray}
q(t)=q_0 \cos(\omega t)+p_0/\omega \sin(\omega t)\\
p(t)=-q_0 \omega \sin(\omega t)+p_0 \cos(\omega t)
\end{eqnarray}
At $T=2\pi /\omega +2k\pi$, we have $q(T)=q_0$, that is $(q,q_0)$ are not independent variables and
the generating function $F_1$ is undefined at this instant. We say that $F_1$ is singular at $T$.
We now prove that at least one of the generating functions is not singular at every instant.

\begin{prop}
Consider the flow $\Phi_t$ of an Hamiltonian system $\phi:(q_0,p_0,t)\mapsto (\Phi_t(q_0,p_0),t)$,
where
$$\Phi_t: (q_0,p_0)\mapsto (\Phi^1_t(q_0,p_0)=q(q_0,p_0,t),\Phi^2_t(q_0,p_0)=p(q_0,p_0,t))$$ For
every $t$, there exists two subsets of cardinal $n$ of the set $(1,\cdots,2n)$, $I_n$ and $K_n$,
such that
\begin{equation}
\det\left(\pfrac{\tilde{\Phi}_{t_{i}}}{z_{j}}\right)_{i\in I_n, j\in K_n}\neq 0
\end{equation}
where $\tilde{\Phi}_t(q,p,q_0,p_0)=(\tilde{\Phi}^1_t(q,q_0,p_0)=q-\Phi^1_t(q_0,p_0),\tilde{\Phi}
^2_t(p,q_0,p_0)=p-\Phi^2_t(q_0,p_0))$ and $z=(q_0,p_0)$
\end{prop}
\begin{proof}
To prove this property, we only need to notice that $\Phi_t$ is a diffeomorphism, i.e.,
$\tilde{\Phi}_{t_{i\in I_n}}$ is an injection, therefore, there exists at least one $n$-dimensional
subspace on which the restriction of $\tilde{\Phi}_{t_{i\in I_n}}$ is a diffeomorphism.
\end{proof}

\begin{theo}
At every instant, at least one generating function is well-defined. Moreover, when they exist,
generating functions define local $C^\infty$-diffeomorphism. \label{theo:existence}
\end{theo}
\begin{proof}
From the previous theorem, there exists $I_n$ and $K_n$ such that
$\det\left(\pfrac{\tilde{\Phi}_{t_i}}{z_j}\right)_{i\in I_n, j\in K_n} \neq 0$. Without any loss of
generality and for simplicity, let the partition be $I_n=(1,\cdots,n)=K_n$,
then we have:
\begin{eqnarray}
\det\left(\pfrac{\tilde{\Phi}^1_t(q,q_0,p_0)}{q_0}\right)&\ne&0
\end{eqnarray}
Moreover, $\tilde{\Phi}^1_t$ verifies:
\begin{equation}\tilde{\Phi}^1_t(q,q_0,p_0)=0\end{equation}
From the local inversion theorem there exists a local diffeomorphism $f_1$ in a neighborhood of
$(q_0,p_0)$ such that $q_0=f_1(q,p_0)$. In addition, the flow defines $p$ as a function of
$(q_0,p_0)$, i.e., replacing $q_0$ by $f_1(q,p_0)$, we obtain $(q_0,p)=(f_1(q,p_0),f_2(q,p_0))$
where $f_2(q,p_0)=\Phi_t^2(f_1(q,p_0),p_0)$. This equation is equivalent to the two equations
verified by $F_2$, hence $f_1=\pfrac{F_2}{p_0}$ and $f_2=\pfrac{F_2}{q}$. This proves that $F_2$
exists and since $\Phi_t$ is $C^\infty$, $F_2$ defines a local $C^\infty$-diffeomorphism from
$(q,p_0)$ to $(p,q_0)$.
\end{proof}
\begin{rem}
The theorem above can be stated for generating functions associated with an arbitrary canonical
transformation, not only the one induced by the phase flow. To proceed the above proof we only
required that the flow defines a $C^\infty$-diffeomorphism, this property is shared by all
canonical transformations.
\end{rem}

Through the harmonic oscillator example, we saw that a generating function may become singular. We
now characterize singularities and give a physical interpretation to them.
\begin{prop}
The generating function $F_{I_p,K_r}$ is singular at time $t$ if and only if
\begin{equation}
\det\left(\pfrac{\tilde{\Phi}_{t_i}}{z_j}\right)_{i \in I, j\in J}= 0
\end{equation}
where $I=\{i\in I_p\}\bigcup \{n+i, i\in \bar{I}_p\}$ and $J=\{j\in \bar{K}_r\}\bigcup \{n+j, j\in
K_r\}$.
\end{prop}
\begin{proof}
The proof proceeds as the previous one, it is also based on local inversion theorem.
\end{proof}

From the above theorem, we deduce that a generating function is singular when there exists multiple
solutions to the boundary value problem. In the harmonic oscillator example, whatever the initial
momentum is, the initial position and position at time $T=2\pi/\omega +2k\pi$ are equal.

Finally, if the Hamiltonian function is independent of time, the system is reversible and therefore
the generating functions $F_{I_p,K_r}$ and $F_{K_r,I_p}$ are similar in the sense that there exists
a diffeomorphism which transforms one into the other. In particular, they develop singularities at
the same instant. If $p=n$ and $r=0$, we obtain that $F_2$ and $F_3$ are similar.

\subsubsection{Linear systems theory}
In this section we particularize the theory developed above to linear systems. The following
developments have implication in the study of relative motion and in optimal control theory as we
will see later. Further, using linear systems theory, we are able to characterize singularities of
generating functions using the state transition matrix.

\paragraph{Hamilton-Jacobi equation}
When studying the relative motion of two particles, one often linearizes the dynamics about the
trajectory of one of the particles (called the nominal trajectory) and then uses a linear
approximation of the dynamics to study the motion of the other particle relative to the nominal
trajectory (perturbed trajectory). Thus, the study of relative motion reduces to the study of a
time-dependent linear Hamiltonian system, i.e., a system with a quadratic Hamiltonian function
without any linear terms\cite{gui02}:
\begin{equation}
H_h=\undemi X_h^T
    \begin{pmatrix}
    H_{qq}(t) & H_{qp}(t)\\
    H_{pq}(t) & H_{pp}(t)
    \end{pmatrix} X_h
\end{equation}
where $X_h=\left(\begin{smallmatrix} \hq\\ \hp \end{smallmatrix}\right)$ is the relative state
vector. Guibout and Scheeres\cite{gui02} proved that the generating functions for the phase flow
transformation must then be quadratic without any linear terms, that is, if we take $F_2$ for
example:
\begin{equation}
F_2=\undemi Y^T
    \begin{pmatrix}
    F^2_{11}(t) & F^2_{12}(t)\\
    F^2_{21}(t) & F^2_{22}(t)
    \end{pmatrix} Y
\label{eq:f2_order2}
\end{equation}
where $Y=\left(\begin{smallmatrix} \hq\\ \kpz
\end{smallmatrix}\right)$ and $\left(\begin{smallmatrix} \kqz\\ \kpz \end{smallmatrix}\right)$
is the relative state vector at initial time. We also point out that both matrices defining $H_h$
and $F_2$ are symmetric. Then Eq.
\ref{eq:f2_in_generala} reads:
\begin{eqnarray}
\hp&=&\pfrac{F_2}{\hq}\nonumber\\
&=&\begin{pmatrix}  F^2_{11}(t)&  F^2_{12}(t)\end{pmatrix}Y
\end{eqnarray}
Substituting into Eq. \ref{eq:f2_in_generalc} yields\footnote{For the canonical transformation
induced by the phase flow, we have seen that $K=0$}:
\begin{equation}
0=Y^T\left\{
\begin{pmatrix}
    \dot{F}^2_{11}(t) & \dot{F}^2_{12}(t)\\
    \dot{F}^2_{12}(t)^T & \dot{F}^2_{22}(t)
\end{pmatrix}
+
\begin{pmatrix}
I& F^2_{11}(t)^T\\
0& F^2_{12}(t)^T
\end{pmatrix}
\begin{pmatrix}
    H_{qq}(t) & H_{qp}(t)\\
    H_{pq}(t)& H_{pp}(t)
\end{pmatrix}
\begin{pmatrix}
I&0 \\
F^2_{11}(t)&F^2_{12}(t)
\end{pmatrix}
\right\}Y
\label{eq:mateq}
\end{equation}
Though the above equations have been derived using $F_2$, they are also valid for $F_1$ (replacing
$Y=\left(\begin{smallmatrix} \hq\\ \kpz
\end{smallmatrix}\right)$ by $Y=\left(\begin{smallmatrix} \hq\\ \kqz
\end{smallmatrix}\right)$) since $F_1$ and $F_2$ solve the same Hamilton-Jacobi
equation (Eqns. \ref{eq:f1_in_generalc}, \ref{eq:f2_in_generalc}). Equation
\ref{eq:mateq} is equivalent to the following $4$ matrix equations:
\begin{equation}
\begin{array}{c}
\dot{F}^{1,2}_{11}(t)+ H_{qq}(t)+H_{qp}(t)F^{1,2}_{11}(t)+F^{1,2}_{11}(t)H_{pq}(t)
+F^{1,2}_{11}(t)H_{pp}(t)F^{1,2}_{11}(t) =0 \\
\\
\dot{F}^{1,2}_{12}(t)+ H_{qp}(t)F^{1,2}_{12}(t)+F^{1,2}_{11}(t)H_{pp}(t)F^{1,2}_{12}(t) =0 \\
\\
\dot{F}^{1,2}_{21}(t)+F^{1,2}_{21}(t)H_{pq}(t)+F^{1,2}_{21}(t)H_{pp}(t)F^{1,2}_{11}(t)  =0\\
\\
\dot{F}^{1,2}_{22}(t)+F^{1,2}_{21}(t)H_{pp}(t)F^{1,2}_{12}(t)=0
\end{array}
\label{eq:matrix_eq_F2}
\end{equation}
where we replaced $F^{2}_{ij}$ by $F^{1,2}_{ij}$ to signify that these equations are valid for both
$F_1$ and $F_2$ and recall that $F^{1,2}_{21}={F^{1,2}_{12}}^T$. A similar set of equations can be
derived for any generating function $F_{I_p,K_r}$, here we only give the equations verified by
$F_3$ and $F_4$:
\begin{equation}
\begin{array}{c}
\dot{F}^{3,4}_{11}(t)+ H_{pp}(t)-H_{pq}(t)F^{3,4}_{11}(t)-F^{3,4}_{11}(t)H_{qp}(t)+
F^{3,4}_{11}(t)H_{qq}(t)F^{3,4}_{11}(t) =0\\
\\
\dot{F}^{3,4}_{12}(t) -H_{pq}(t)F^{3,4}_{12}(t)+F^{3,4}_{11}(t)H_{qq}(t)F^{3,4}_{12}(t) =0   \\
\\
\dot{F}^{3,4}_{21}(t)-F^{3,4}_{21}(t)H_{qp}(t)+F^{3,4}_{21}(t)H_{qq}(t)F^{3,4}_{11}(t)  =0\\
\\
\dot{F}^{3,4}_{22}(t)+ F^{3,4}_{21}(t)(t)H_{qq}(t)F^{3,4}_{12}(t)=0\end{array}
\label{eq:matrix_eq_F3}
\end{equation}

The first equations of Eqns \ref{eq:matrix_eq_F2} and \ref{eq:matrix_eq_F3} are Ricatti equations,
the second and third are non-homogeneous, time varying, linear equations and are equivalent to each
other (i.e., transform into each other under transpose), and the last are just a quadrature.

\paragraph{Perturbation matrices}
Another approach to the study of relative motion of spacecraft is to use the state transition
matrix. This method was developed by Battin \cite{bat99} for the case of a spacecraft moving in a
point mass gravity field. Let $\Phi$ be the state transition matrix which describes the relative
motion:
\begin{equation}
\left(\begin{array}{c} \hq\\ \hp \end{array}
\right)=\Phi\left(\begin{array}{c} \kqz\\ \kpz \end{array} \right)
\end{equation}
where $\Phi=\left(\begin{array}{cc} \Phi_{qq}&\Phi_{qp}\\
\Phi_{pq}&\Phi_{pp} \end{array} \right)$.

From the state transition matrix, Battin\cite{bat99} defines the fundamental perturbation matrices
$C$ and $\tilde{C}$ as:
\begin{eqnarray}
\tilde{C}&=&\Phi_{pq}\Phi_{qq}^{-1}\nonumber\\
C&=&\Phi_{pp}\Phi_{qp}^{-1}
\end{eqnarray}
That is, given $\kpz=0$, $\tilde{C}\hq=\hp$ and given $\kqz=0$, $C\hq=\hp$. He shows that for
relative motion of a spacecraft in a point mass gravity field these matrices verify a Ricatti
equation and are therefore symmetric. Using the generating functions for the canonical
transformation induced by the phase flow, we immediately recover these properties and also show
that they are verified for any relative motion of two particles in a Hamiltonian dynamical system.

From Eqns. \ref{eq:f2_in_generala} and \ref{eq:f2_in_generalb}:
\begin{eqnarray}
\hp&=&\pfrac{F_2}{\kpz}\nonumber\\
   &=&F^2_{11} \hq+F^2_{12} \kpz\\
\kqz&=&\pfrac{F_2}{\hq}\nonumber\\
     &=&F^2_{21} \hq+F^2_{22} \kpz
\end{eqnarray}

Solving for $(\hq, \hp)$ yields:
\begin{align}
&\hq={F^2_{21}}^{-1}\kqz-{F^2_{21}}^{-1}F^2_{22}\kpz\\
&\hp=F^2_{11}{F^2_{21}}^{-1}\kqz+(F^2_{12}-F^2_{11}{F^2_{21}}^{-1}F^2_{22})\kpz
\end{align}
From the above equations we are able to link the state transition matrix to the generating function
$F_2$.

\begin{displaymath}
\left\{
\begin{array}{ccl}
\Phi_{qp}&=& -{F^2_{21}}^{-1}F^2_{22}\\
\Phi_{qq}&=&  {F^2_{21}}^{-1}\\
\Phi_{pp}&=& F^2_{12}-F^2_{11}{F^2_{21}}^{-1}F^2_{22}\\
\Phi_{pq}&=& F^2_{11}{F^2_{21}}^{-1}
\end{array}
\right.
\end{displaymath}
We conclude that
\begin{equation}
\tilde{C}=\Phi_{pq}\Phi_{qq}^{-1}=F^2_{11}
\end{equation}

In the same manner, but using $F_1$, we can show that:
\begin{equation}
C=\Phi_{pp}\Phi_{qp}^{-1}=F^1_{11}
\end{equation}

Thus, we have shown that $C$ and $\tilde{C}$ are symmetric by nature (as $F^{1,2}_{11}$ is
symmetric by definition),
 and moreover that they verify the Ricatti equation given in Eq.
\ref{eq:matrix_eq_F2}.

\paragraph{Singularities of generating functions and state transition matrix}
From Eqns. \ref{eq:f2_in_generala}, \ref{eq:f2_in_generalb} and \ref{eq:f2_order2}, we derive a
relationship between terms of $F_2$ and some coefficients of the state transition matrix:
\begin{eqnarray}
\hp&=&\pfrac{F_2}{\hp}\nonumber\\
&=&F^2_{11} \hq +F^2_{12} \kpz \nonumber\\
&&\textrm{but we also have}\nonumber\\
\hp&=&\Phi_{pq}\Phi_{qq}^{-1}\hq
+(\Phi_{pp}-\Phi_{pq}\Phi_{qq}^{-1}\Phi_{qp})\kpz \nonumber\\
&&\\
\kqz&=&\pfrac{F_2}{\kpz} \nonumber\\
&=&F^2_{21} \hq + F^2_{22} \kpz \nonumber\\
&&\textrm{but we also have}\nonumber\\
\kqz&=&\Phi_{qq}^{-1}\hq-\Phi_{qq}^{-1}\Phi_{qp}\kpz
\end{eqnarray}
Thus:
\begin{eqnarray}
F^2_{11}&=&\Phi_{pq}\Phi_{qq}^{-1}\\
F^2_{12}&=&\Phi_{pp}-\Phi_{pq}\Phi_{qq}^{-1}\Phi_{qp}\\
F^2_{21}&=&\Phi_{qq}^{-1}\\
F^2_{22}&=&\Phi_{qq}^{-1}\Phi_{qp}\\
\end{eqnarray}

We conclude that if $\Phi_{qq}$ is singular, $F_2$ is also singular. The same analysis can be
achieved for the other generating functions, and we find that:
\begin{itemize}
\item $F_1$ is singular when $\Phi_{qp}$ is singular,
\item $F_3$ is singular when $\Phi_{pp}$ is singular,
\item $F_4$ is singular when $\Phi_{pq}$ is singular.
\end{itemize}
These results can be extended to other generating functions, but requires us to work with another
block decomposition of the state transition matrix.

\subsubsection{On singularities of generating functions}
\label{sect:singularities} We have proved local existence of generating functions and mentioned
that they may not be globally defined for all time. Using linear system theory we were able to
predict where the singularities are and to interpret their meaning as multiple solutions to the
two-point boundary value problem. In this section we extend our study to singularities of nonlinear
systems.

\paragraph{Lagrangian submanifold} Consider an arbitrary generating function $F_{I_p,K_r}$. Then Eqns.
\ref{eq:fjpr_in_generala}-\ref{eq:fjpr_in_generald}
define a $2n$-dimensional submanifold called a canonical relation\cite{wei77} of the
$4n$-dimensional symplectic space $P_1 \times P_2$. In addition, since the new variables $(Q,P)$
(or $(q_0,p_0)$) do not appear in the Hamilton-Jacobi equation
\ref{eq:fjpr_in_general_hj} we may consider them as parameters. In that case Eqns.
\ref{eq:fjpr_in_generala} and \ref{eq:fjpr_in_generalb} define an $n$-dimensional submanifold of the
symplectic space $P_1$ called a Lagrangian submanifold\cite{wei77}. The study of singularities can
be achieved using either canonical relations\cite{abr78} or Lagrangian
submanifolds\cite{arn88a,mar98}.
\begin{theo}
The generating function $F_{I_p,K_r}$ is singular if and only if the local projection of the
canonical relation $\cal L$ defined by Eqns. \ref{eq:fjpr_in_generala}-\ref{eq:fjpr_in_generald}
onto $(q_{I_p},p_{\bar{I}_p},Q_{K_r},P_{\bar{K}_r})$ is not a local diffeomorphism.
\end{theo}
Moreover, the projection of such a singular point is called a caustic. If one works with Lagrangian
submanifolds then the previous theorem becomes
\begin{theo}
The generating function\footnote{We consider here that the generating function is function of $n$
variables only, and has $n$ parameters.} $F_{I_p,K_r}$ is singular if the local projection of the
Lagrangian submanifold defined by Eqns.
\ref{eq:fjpr_in_generala} and
\ref{eq:fjpr_in_generalb} onto $(q_{I_p},p_{\bar{I}_p})$ is not a local diffeomorphism.
\end{theo}

In light of these previous theorems, we can give a geometrical interpretation to theorem
\ref{theo:existence} on the existence of generating functions. Given a canonical relation $\cal L$
(or a Lagrangian submanifold) defined by a canonical transformation, there exists a
$2n$-dimensional (or $n$-dimensional) submanifold $\cal M$ of $P_1 \times P_2$ (or $P_1$) such that
the local projection of $\cal L$ onto $\cal M$ is a local diffeomorphism.

\paragraph{Study of caustics}
To study caustics two approaches, at least, are possible depending on the problem. A good
understanding of the physics may provide information very easily. For instance, consider the two
body problem in dimension $2$, and the problem of going from one point $A$ to a point $B$,
symmetric with respect to the central body, in a certain lapse of time. The trajectory that links
$A$ to $B$ is an ellipse whose perigee and apogee are $A$ and $B$. Therefore, there are two
solutions to this problem depending upon which way the particle is going. In terms of generating
functions, we deduce that $F_3$ is nonsingular (there is a unique solution once the final momentum
is given) but $F_1$ is singular (existence of two solutions) and the caustic is a
fold\footnote{Maps from $\mathbb{R}^2$ into $\mathbb{R}^2$ have two types of stable singularities,
folds and cusps. However, only folds have two antecedents, cusps have three.}. The other method to
study caustics consists in using a known nonsingular generating function to define the Lagrangian
submanifold $\cal L$ and then study its projection. A very illustrative example is given by Ehlers
and Newman
\cite{ehl00new}, they treat the evolution of an ensemble of free particles whose initial momentum
distribution is $p=\frac{1}{1+q^2}$ using the Hamilton-Jacobi equation and generating functions for
the phase flow canonical transformation. They are able to solve the problem analytically, that is,
identify a time at which $F_1$ is singular, find the equations defining the Lagrangian submanifold
using $F_3$ and study its projection to eventually find two folds. Nevertheless, such an analysis
is not always possible as solutions to the Hamilton-Jacobi equation are usually found numerically,
not analytically. In the remainder of this section, we focus on a class of problem that can be
solved numerically for which we are able to characterize the caustics.

Suppose we are interested in the relative motion of a particle, called the deputy, whose
coordinates are $(q,p)$ with respect to another one, called the chief, whose coordinates are
$(q^0,p^0)$, both moving in an Hamiltonian field. If both particles stay ``close" to each other, we
can expand $(q,p)$ as a Taylor series about the trajectory of the chief. The dynamics of the
relative motion is described by the Hamiltonian function $H_h$\cite{gui02}:
\begin{equation}
H^h(X^h,t)=\sum_{p=2}^\infty \sum_{\begin{array}{c}i_1,\cdots,i_{2n}=0 \\
i_1+\cdots+i_{2n}=p
\end{array}}^p
\frac{1}{i_1!\cdots i_{2n}!}
\frac{\partial^{p} H}{\partial q_1^{i_1}\cdots\partial q_n^{i_n}
\partial p_1^{i_{n+1}}\cdots\partial p_n^{i_{2n}}}(q^0,p^0,t) {X^h_1}^{i_1}\dots{X^h_{2n}}^{i_{2n}}
\label{eq:l_hamilton_n}
\end{equation}
where  $X^h=(\hq,\hp)$,  $\hq=q-q^0$ and $\hp=p-p^0$. Since $H_h$ has infinitely many terms, we are
usually not able to solve the Hamilton-Jacobi equation but we can approximate the dynamics by
truncating the series $H_h$ in order to only keep finitely many terms. Suppose $N$ terms are kept,
then we say that we describe the relative motion using an approximation of order $N$. Clearly, the
greater $N$ is, the better our approximation is to the nonlinear motion of a particle about the
nominal trajectory. When an approximation of order $N$ is used, we look for a generating function
$F_{I_p,K_r}$ as a polynomial of order $N$ in its spatial variables with time dependent
coefficients. The Hamilton-Jacobi equation reduces to a set of ordinary differential equations that
we integrate numerically. Once $F_{I_p,K_r}$ is known, we find the other generating functions from
the Legrendre transformation, at the cost of a series inversion. If a generating function is
singular, the inversion does not have a unique solution, the number of solutions characterizes the
caustic. To illustrate this method, let us consider the following example.

\paragraph{Motion about the Libration point $L_2$ in the Hill three-body problem} Consider a
spacecraft moving about and staying close to the Libration point $L_2$ in the Hill three-body
problem (See the appendix for a description of the Hill three-body problem). Its relative motion
with respect to $L_2$ is described by the Hamiltonian function $H_h$ (Eq. \ref{eq:l_hamilton_n})
and approximated at order $N$ by truncation of terms of order greater than $N$ in the Taylor series
defining $H_h$. Using the algorithm developed by Guibout and Scheeres\cite{gui02} we find the
generating functions for the canonical transformation induced by the approximation of the phase
flow, that is, the Taylor series expansion up to order $N$ of the exact generating function about
the Libration point $L_2$.
\begin{eqnarray}
F_2(q_x,q_y,p_{0_x},p_{0_y},t)&=&
f^2_{11}(t) q_x^2 +f^2_{12}(t) q_x q_y+f^2_{13}(t) q_x p_{0_x}+f^2_{14}(t) q_x p_{0_y}\nonumber\\
&&f^2_{22}(t) q_y^2+f^2_{23}(t) q_y p_{0_x}(t)+f^2_{24}(t) q_y p_{0_y} \nonumber\\
&&f^2_{33}(t) p_{0_x}^2+f^2_{34}(t) p_{0_x} p_{0_y}+f^2_{44}(t)
p_{0_y}^2+r(q_x,q_y,p_{0_x},p_{0_y},t)
\end{eqnarray}
where $(q,p,q_{0},p_{0})$ are relative position and momenta of the spacecraft with respect to $L_2$
at $t$ and at $t_0$, the initial time, and $r$ is a polynomial of degree $N$ in its spatial
variables with time dependent coefficients and without any quadratic terms. At $T=1.6822$, $F_1$ is
singular but $F_2$ is not. Eqns.
\ref{eq:f2_in_generala} and \ref{eq:f2_in_generalb} reads:
\begin{eqnarray}
p_x&=& 2 f^2_{11}(T) q_x+f^2_{12}(T) q_y+f^2_{13}(T) p_{0_x}+f^2_{14}(T) p_{0_y}+
D_1 r(q_x,q_y,p_{0_x},p_{0_y},T)\label{eq:f21}\\
p_y&=&f^2_{12}(T) q_x +2 f^2_{22}(T) q_y+f^2_{23}(T) p_{0_x}+f^2_{24}(T) p_{0_y}+
D_2 r(q_x,q_y,p_{0_x},p_{0_y},T)\label{eq:f22}\\
q_{0_x}&=& f^2_{13}(T) q_x + f^2_{23}(T) q_y + 2 f^2_{33}(T) p_{0_x} + f^2_{34}(T) p_{0_y}+
D_3 r(q_x,q_y,p_{0_x},p_{0_y},T)\label{eq:f23}\\
q_{0_y}&=& f^2_{14}(T) q_x + f^2_{24}(T) q_y + f^2_{34}(T) p_{0_x} + 2 f^2_{44}(T) p_{0_y}+ D_4
r(q_x,q_y,p_{0_x},p_{0_y},T)\label{eq:f24}
\end{eqnarray}
where $D_i r$ represents the derivative of $r$ with respect to its $i^{th}$ variable. Eqns.
\ref{eq:f21}-\ref{eq:f24} define a canonical relation $\cal L$. By assumption $F_1$ is singular,
therefore the projection of $\cal L$ onto $(q,q_0)$ is not a local diffeomorphism and there exists
a caustic. The theory developed above provides a technique to study this caustic using $F_2$. Eqns.
\ref{eq:f21}-\ref{eq:f24} provide $p$ and $q_0$ as a function of $(q,p_0)$, but to characterize the
caustic we need $p$ and $p_0$ as a function of $(q,q_0)$. $F_1$ being singular, there are multiple
solutions to this problem, and one valuable piece of information is the number $k$ of such
solutions. To find $p$ and $p_0$ as a function of $(q,q_0)$ we can first invert equations
\ref{eq:f23} and
\ref{eq:f24} to express $p_0$ as a function of $(q,q_0)$ and then plug this relation into Eqns.
\ref{eq:f21} and \ref{eq:f22}. The first step requires a series inversion that can be proceeded
using the technique developed in \cite{mou30} by Moulton. Let us rewrite Eqns. \ref{eq:f23} and
\ref{eq:f24}:
\begin{eqnarray}
2 f^2_{33}(T) p_{0_x} + f^2_{34}(T) p_{0_y} &=& q_{0_x} - f^2_{13}(T) q_x - f^2_{23}(T) q_y -
D_3 r(q_x,q_y,p_{0_x},p_{0_y},T)\label{eq:f25}\\
f^2_{34}(T) p_{0_x} + 2 f^2_{44}(T) p_{0_y}&=& q_{0_y} - f^2_{14}(T) q_x - f^2_{24}(T) q_y - D_4
r(q_x,q_y,p_{0_x},p_{0_y},T)\label{eq:f26}
\end{eqnarray}
The determinant of the coefficients of the linear terms on the left hand side is zero (otherwise
there is a unique solution to the series inversion) but each of the coefficients is non zero, that
is, we can solve for $p_{0_x}$ as a function of $(p_{0_y},q_{0_x},q_{0_y})$ using equation
\ref{eq:f25}. Then we plug this solution into Eq. \ref{eq:f26} and we obtain an equation of the
form
\begin{equation}
R(p_{0_y},q_{0_x},q_{0_y})=0 \label{eq:f27}
\end{equation}
that contains no terms in $p_{0_y}$ alone of the first degree. In addition, $R$ contains a non zero
term of the form $\alpha p_{0_y}^2$, where $\alpha$ is a real number. In this case, Weierstrass
proved that there exist $2$ solutions $p^1_{0_y}$ and $p^2_{0_y}$ to Eq. \ref{eq:f27}, that is, the
caustic is a fold.

In the same way, we can study the singularity of $F_1$ at initial time. At $t=0$, $F_2$ generates
the identity transformation, hence $f^2_{33}(0)= f^2_{34}(0)=f^2_{34}(0)= f^2_{44}(0)=0$. This time
there is no nonzero first minor, and we find that there exists infinitely many solutions to the
series inversion. Another way to see this is to use the Legendre transformation:
\begin{equation}
F_1(q,q_0,t)=F_2(q,p_0,t)-q_0 p_0
\end{equation}
As $t$ tends toward $0$, $(q,p)$ goes to $(q_0,p_0)$ and $F_2$ converges toward the identity
transformation $q p_0\rightarrow_{t\rightarrow0} q_0 p_0$. Therefore, as $t$ goes to $0$, $F_1$
also goes to $0$, i.e., the projection of $\cal L$ onto $(q,q_0)$ reduces to a point.

There are many other results on caustics and Lagrangian submanifolds that go beyond the scope of
this paper. Study of the Lagrangian submanifold at singularities is ``the beginning of deep
connections between symplectic geometry, geometric optics, partial differential equations, and
Fourier integral operators." (R. Abraham \cite{abr78}), we refer to Abraham \cite{abr78} and
references given therein for more information on this subject. Let us now come back to two-point
boundary value problems.

So far we have studied the generating functions associated with the canonical transformation
induced by the phase flow and showed they formally solve any two-point boundary value problem.
Nonetheless, for Hamiltonian dynamical systems there exists another function, called Hamilton's
principal function, that solves the same problem and thus for completeness we discuss it. In this
section we introduce this function and show how it compares to the generating functions for the
canonical transformation induced by the phase flow.
\subsection{Hamilton's principal function}
Though generating functions are used in this paper to solve boundary value problems, they have been
introduced by Jacobi and mostly used thereafter as fundamental functions which can yield all the
equations of motion by simple differentiations and eliminations, without integration. Nevertheless,
it was Hamilton who first hit upon the idea of finding such a fundamental function, he proved its
existence in geometrical optics (i.e., for time independent Hamiltonian systems) in $1834$ and
called it characteristic function\cite{ham34}. The year later, he published a second
essay\cite{ham35} on systems of attracting and repelling points in which he showed that the
evolution of dynamical systems is characterized by a single function called Hamilton's principal
function: ``The former Essay contained a general method for reducing all the most important
problems of dynamics to the study of one characteristic function, one central or radical relation.
It was remarked at the close of that Essay, that many eliminations required by this method in its
first conception, might be avoided by a general transformation, introducing the time explicitly
into a part S of the whole characteristic function V ; and it is now proposed to fix the attention
chiefly on this part S, and to call it the Principal Function." (William R. Hamilton, in the
introductory remarks of ``Second essay on a General Method in Dynamics" \cite{ham35})

\subsubsection{Hamilton's principal function to describe the phase flow}
As with generating functions, Hamilton's principal function may be derived using the calculus of
variations. Consider the extended action integral:
\begin{equation}
A=\int_{\tau_0}^{\tau_1}( pq' +p_t t' )d\tau \label{eq:action}
\end{equation}
under the auxiliary condition $K(q,t,p,p_t)=0$, where $q'=dq/d\tau$, $p_t$ is the momentum
associated with the generalized coordinates $t$ and $K=p_t+H$.

Define a line element\footnote{Note that the geometry established by this line element is not
Riemannian\cite{lan77}} $d\sigma$ for the extended configuration space $(q,t)$ by
\begin{equation}
d\sigma=Ldt=Lt'd\tau
\end{equation}
Then, we can connect two points $(q_0,t_0)$ and $(q_1,t_1)$ of the extended configuration space by
a shortest line $\gamma$ and measure its length from:
\begin{equation}
A=\int_\gamma d\sigma=\int_\gamma Lt'd\tau
\end{equation}
The distance we obtain is function of the coordinates of the end-points and is called Hamilton's
principal function: $W(q_0,t_0,q_1,t_1)$.

From calculus of variations\cite{lan77} we know that the variation of the action $A$ can be
expressed as a function of the boundary terms if we vary the limits of the integral:
\begin{eqnarray}
\delta A =p_0\delta q_0+p_{t_0} \delta t_0 -p_1\delta q_1-p_{t_1}\delta t_1
\end{eqnarray}
On the other hand we have:
\begin{equation}
\delta A =\delta W(q_0,t_0,q_1,t_1)=\pfrac{W}{q_0}\delta q_0+\pfrac{W}{t_0}\delta t_0
+\pfrac{W}{q_1}\delta q_1+\pfrac{W}{t_1}\delta t_1
\end{equation}
that is:
\begin{eqnarray}
p_0&=&\pfrac{W}{q_0}(q_0,t_0,q_1,t_1)\label{eq:W1}\\
p_1&=&-\pfrac{W}{q_1}(q_0,t_0,q_1,t_1)\label{eq:W2}
\end{eqnarray}
and
\begin{eqnarray}
\pfrac{W}{t_0}(q_0,t_0,q_1,t_1)+H(q_0,\pfrac{W}{q_0},t_0)&=&0\label{eq:pde1}\\
-\pfrac{W}{t_1}(q_0,t_0,q_1,t_1)+H(q_1,-\pfrac{W}{q_1},t_1)&=&0\label{eq:pde2}
\end{eqnarray}
where $K$ has been replaced by $p_t+H$. As with generating functions of the first kind, Hamilton's
principal function solves boundary value problems of Lambert's type through Eqns. \ref{eq:W1} and
\ref{eq:W2}. To find $W$, however, we need to solve a system of two partial differential equations
(Eqns. \ref{eq:pde1} and \ref{eq:pde2}).

\subsubsection{Hamilton's principal function and generating functions}
In this section we highlight the main differences between generating functions for the canonical
transformation induced by the phase flow and Hamilton's principal function. For sake of simplicity
we compare $F_1(q,q_0,t)$ and $W(q,t,q_0,t_0)$.
\paragraph{Calculus of variation} Even if both functions are derived from calculus of variations,
there are fundamental differences between them. To derive generating functions we used the
principle of least action with the time $t$ as independent variables whereas we increase the
dimensionality of the system by adding the time $t$ to the generalized coordinates to derive
Hamilton's principal function. As a consequence, generating functions generates a transformation
between two points in the phase space, i.e., they act without passage of time whereas Hamilton's
principal function generates a transformation between two points in the extended phase space, i.e.,
between two points in the phase space with different times. This difference may be viewed as
follows: Generating functions allow to characterize the phase flow given an initial time, $t_0$
(i.e., to characterize all trajectories whose initial conditions are specified at $t_0$), whereas
Hamilton's principal function does not impose any constraint on the initial time. The counterpart
being that Hamilton's principal function must satisfy two partial differential equations (Eq.
\ref{eq:pde1} defines $W$ as a function of $t_0$ and Eq. \ref{eq:pde2} defines $W$ as a function of
$t_1$) whereas generating functions satisfy only one.

Moreover, to derive the generating functions fixed endpoints are imposed, that is we impose the
trajectory in both sets of variables to verify the principle of least action. On the other hand,
the variation used to derive Hamilton's principal function involves moving endpoints and an energy
constraint. This difference may be interpreted as follows: Hamilton's principal function generates
a transformation which maps a point of a given energy surface to another point on the same energy
surface and is not defined for points that do not lie on this surface. As a consequence of the
energy constraint we have:
\begin{equation}
|\ppfrac{W}{q_0}{q_1}|=0 \label{eq:Wsing}
\end{equation}
As noticed by Lanczos\cite{lan77}, ``this is a characteristic property of the $W$-function which
has no equivalent in Jacobi's theory". On the other hand, generating functions map any point of the
phase space into another one, the only constraint is imposed through the principle of least action
(or equivalently by the definition of canonical transformation): we impose the trajectory in both
sets of coordinates to be Hamiltonian with Hamiltonian function $H$ and $K$ respectively.

\paragraph{Fixed initial time} In the derivation of Hamilton's principal function $dt_0$ may be
chosen to be zero, that is, the initial time is imposed. Hamilton's principal function loses its
dependence with respect to $t_0$, Eq. \ref{eq:pde1} is trivially verified and Eq. \ref{eq:Wsing}
does not hold anymore, $W$ and $F_1$ become equivalent.

Finally, in \cite{ham35} Hamilton also derives another principal function $Q(p_0,t_0,p_1,t_1)$
which compares to $W$ as $F_4$ compares to $F_1$, the derivation being the same we will not go
through it.

To conclude, Hamilton's principal function appears to be more general than the generating functions
for the canonical transformation induced by the phase flow. On the other hand, to solve a two-point
boundary value problem, initial and final times are specified and therefore, any of these functions
will identically solve the problem. To find Hamilton's principal function, we need to solve two
partial differential equations whereas only one need to be solved to find the generating functions.
For this reason, generating functions will be used in the following examples.

\section{Applications}
\subsection{Solving the optimal control problem using the generating functions}
The use of the generating functions to solve an optimal control problem has first been addressed by
Scheeres, Guibout and Park\cite{sch03gui}. They suggested an indirect approach for evaluating the
initial adjoints without initial guess. In the present paper,  we review their approach and
generalize it to a wider class of problem.

\paragraph{Problem formulation}
Assume a dynamical system described by:
\begin{eqnarray}
\dot{x}&=&f(x,u,t)
\label{eq:dynamics}\\
x(t=0)&=&x_0
\label{eq:initialcond}
\end{eqnarray}
where $u$ is the control variable, $x\in \mathbb{R}^n$ and $u\in\mathbb{R}^m$. An optimal control
problem is formulated as follows:
\begin{equation}
\min_u K(x(t_f))+\int_{t_0}^{t_f} L(x,u,t) dt
\label{eq:optimal}
\end{equation}
where $t_f$ is the known final time. This formulation is called the Bolza formulation. Other
formulations are possible and completely equivalent
\begin{eqnarray}
\min_u \tilde{K}(x(t_f))& \qquad & \textrm{Mayer formulation} \label{eq:mayer}\\
\min_u \int_{t_0}^{t_f} \tilde{L}(x,u,t) dt && \textrm{Lagrange formulation}
\label{eq:lagrange}
\end{eqnarray}

Further, some final conditions may be specified. For instance, suppose that $k$ final conditions
are given for the final state, i.e.,
\begin{equation}
\psi_j(x(t_f),t_f)=0 \quad j\in (1\cdots k)
\label{eq:opt_control_fc}
\end{equation}

\paragraph{Necessary conditions}
Define the Hamiltonian function $H$:
\begin{equation}
H(x,p,u,t)=p^T \dot{x}+L(x,u,t)
\label{eq:h}
\end{equation}
where $p \in \mathbb{R}^n$ is the costate vector. Applying the Pontryagin principle allows one to
find the optimal control:
\begin{equation}
\bar{u}=\arg \min_u H(x,p,u,t)
\label{eq:neccond3}
\end{equation}
Then the necessary conditions for optimality are given by:
\begin{eqnarray}
\dot{x}&=&\pfrac{H}{p}(x,p,\bar{u},t)\label{eq:neccond1}\\
\dot{p}&=&-\pfrac{H}{x}(x,p,\bar{u},t)\label{eq:neccond2}
\end{eqnarray}
To integrate these $2n$ differential equations we need $2n$ boundary conditions: $n+k$ are
specified in the problem statement, the other $n-k$ are given by the transversality conditions:
\begin{equation}
p(t_f)-\pfrac{K}{x}(t_f)= \nu^T \pfrac{\psi}{x}(x(t_f))
\label{eq:transversality}
\end{equation}
where $\nu$ is a $k$-dimensional vector.

\paragraph{Solving the optimal control using the generating functions}
In the following, we are making two assumptions which may be relaxed in future research.
\begin{enumerate}
\item One can solve for $u$ as a function of $(x,p)$ using Eq. \ref{eq:neccond3}, that is, we can
define a new Hamiltonian function $\bar{H}(x,p,t)=H(x,p,\bar{u}(x,p,t),t)$.
\item One can eliminate the $\nu$'s in Eq. \ref{eq:transversality}, so that Eq.
\ref{eq:transversality} becomes
\begin{equation}
p_i(t_f)=p_{f_i} \qquad \forall i \in (k,n)
\end{equation}
and transform Eq. \ref{eq:opt_control_fc} into:
\begin{equation}
x_j(t_f))=x_{f_j} \quad j\in (1\cdots k)
\end{equation}
\end{enumerate}
Then, solving the optimal control problem is equivalent to find the solutions $(x,p)$ satisfying:
\begin{eqnarray}
\dot{x}&=&\pfrac{\bar{H}}{p}(x,p,t) \\
\dot{p}&=&-\pfrac{\bar{H}}{x}(x,p,t)
\end{eqnarray}
with boundary conditions
\begin{equation}
\begin{array}{ccc}
x(t=0)&=&x_0 \quad \quad \quad \quad \qquad \qquad \quad \\
x_i(t_f)&=&x_{f_i} \quad \forall i\in (1,\cdots, k) \\
p_i(t_f))&=&p_{f_i} \quad \forall i \in (k,\cdots,n)\end{array}
\label{eq:boundarycond}
\end{equation}

These equations define a two-point boundary value problem and hence are usually difficult to solve
because they generally require an estimate of the initial (or final) state, which usually has no
physical interpretation. An indirect approach can be developed to solve this problem, namely, the
use of the generating function
$F_{I_n,K_k}(x_{0_1},\cdots,x_{0_n},x_{f_1},\cdots,x_{f_k},p_{f_{k+1}},\cdots,p_{f_n})$. Eqns.
\ref{eq:fjpr_in_generala}, \ref{eq:fjpr_in_generalb} and \ref{eq:fjpr_in_generalc} solves the
boundary value problem and hence the optimal control problem:
\begin{eqnarray}
p_{0_i}&=&-\pfrac{F_{I_n,K_k}}{x_{0_i}} \\
x_{f_i}&=&-\pfrac{F_{I_n,K_k}}{p_{f_i}} \\
p_{f_i}&=&\pfrac{F_{I_n,K_k}}{x_{f_i}}
\end{eqnarray}

In the case where $k=n$, that is initial and final states of the system are specified, the
generating function that must be used to solve the boundary value problem is $F_1$. In that case,
Park and Scheeres\cite{cha03sch} showed that $F_1$ satisfies the Hamilton-Jacobi-Bellmann equation.

\paragraph{Particular case: The linear quadratic problem}
Assume the dynamics of the system is linear:
\begin{equation}
\dot{x}(t)=A(t)x(t)+B(t)u(t)
\end{equation}
and the cost function $J$ is quadratic:
\begin{equation}
J=\frac{1}{2}[Mx(t_f)-m_f]^T Q_f[Mx(t_f)-m_f]+\frac{1}{2}\int_{t_0}^{t_f} \begin{pmatrix} x^T
&u^T\end{pmatrix}\begin{pmatrix}Q &N\\N^T&R\end{pmatrix}\begin{pmatrix} x
\\u\end{pmatrix}
\end{equation}
and $Q$ is symmetric positive semi-definite, $R$ and $Q_f$ are symmetric positive definite.
Moreover, define $L$ to be $L=\frac{1}{2}\begin{pmatrix} x^T &u^T\end{pmatrix}\begin{pmatrix}Q
&N\\N^T&R\end{pmatrix}\begin{pmatrix} x
\\u\end{pmatrix}$

Using previous notations, we define the Hamiltonian function $H$:
\begin{equation}
H(x,p,u)=p^T \dot{x}+L(x,u)
\end{equation}
From equation \ref{eq:neccond3}, we get
\begin{equation}
\bar{u}=-R^{-1}B^Tp-R^{-1}N^Tx \label{eq:u}
\end{equation}
Substituting $\bar{u}$ in Eqns. \ref{eq:neccond1} and \ref{eq:neccond2} yields:
\begin{equation}
\bar{H}(x,p)=H(x,p,-R^{-1}B^Tp-R^{-1}N^Tx)
\end{equation}
and
\begin{eqnarray}
\dot{x}&=&Ax+B(-R^{-1}B^Tp-R^{-1}N^Tx) \label{eq:neccond1_linear}\\
\dot{p}&=&-(A^Tp+Qx+N(-R^{-T}B^Tp-R^{-1}N^Tx))\label{eq:neccond2_linear}
\end{eqnarray}
Boundary conditions for this problem are still given by equations \ref{eq:boundarycond}. Since the
Hamiltonian function defining this system is quadratic, this problem is often solved using the
state transition matrix. We have seen previously that, in linear systems theory, both generating
functions and the state transition matrix are equivalent. Moreover, to compute the generating
function or the state transition matrix, four matrix equations of dimension $n$ need to be solved.
Therefore, both methods are exactly equivalent for the linear quadratic problem. Finally, another
method to solve the linear quadratic problem is to apply Ricatti transformation to reduce the
problem to two matrix ordinary differential equations, a Ricatti equation and a time-varying linear
equation. An analogy can be drawn between these two equations and the ones verified by the
generating function.

\subsection{Finding periodic orbits using the generating functions}
Another application of the generating functions for the canonical transformation induced by the
phase flow is to search for periodic orbits. This application was first presented by Guibout and
Scheeres\cite{gui03b}, we review their methodology in this paper and refer to\cite{gui03b} for more
details and additional examples.

\subsubsection{The search for periodic orbits: a two-point boundary value problem}

The main idea is to transform the search for periodic orbits into a two-point boundary value
problem that can be handled using generating functions. For a periodic orbit of period $T$, both
position and momentum take the same values at $t$ and at $t+kT,\quad k\in\mathbb{Z}$. In terms of
initial conditions, this reads:
\begin{eqnarray}
q(T)&=&q_0\label{eq:bound_eqa}\\
p(T)&=&p_0
\label{eq:bound_eqb}
\end{eqnarray}

For a dynamical system with $n$ degrees of freedom Eqns. \ref{eq:bound_eqa} and \ref{eq:bound_eqb}
can be viewed as $2n$ equations of $2n+1$ variables, the initial conditions $(q_0,p_0)$ and the
period $T$. To solve such a problem, for each trial $(q_0,p_0,T)$ one needs to integrate the
equations of motion and check if the $2n$ equations are verified, and if they are not try again. On
the other hand, Eqns. \ref{eq:bound_eqa} and \ref{eq:bound_eqb}  can also be viewed as a two-point
boundary value problem. Suppose the initial momentum $p_0$ and the position at time $T$, $q$, are
given, then Eqns. \ref{eq:bound_eqa} and \ref{eq:bound_eqb} define $2n$ equations with $2n+1$
variables, the initial position $q_0$, the momentum at time $T$, $p$, and the time period $T$.
Solutions to these equations characterize all periodic orbits. The idea now is to use the
generating functions for the phase flow transformation to solve this problem. Depending on the
two-point boundary value problem we choose to characterize periodic orbits, different generating
functions can be used. In the following we will only deal with generating functions of the first
and second kind, but this theory can be readily generalized to any kind of generating functions.

\subsubsection{Solving the two-point boundary value problem}
\paragraph{Generating functions of the first kind}
The generating function $F_1$ allows us to solve a two-point boundary value problem for which
initial position and position at time $T$ are given. The solution to this problem is found using
Eqns. \ref{eq:f1_in_generala} and
\ref{eq:f1_in_generalb}.
\begin{eqnarray}
p&=&\pfrac{F_1}{q}(q,q_0,T)\\
p_0&=&-\pfrac{F_1}{q_0}(q,q_0,T)
\end{eqnarray}

On the other hand, the boundary value problem that characterizes periodic orbits is defined by
equations
\ref{eq:bound_eqa} and
\ref{eq:bound_eqb}. Hence, combining these four equations yields:
\begin{eqnarray}
  p_0&=&-\pfrac{F_1}{q_0}(q=q_0,q_0,T)\label{eq:per_f1_suffa}\\
p(T)&=&\pfrac{F_1}{q}(q=q_0,q_0,T)\label{eq:per_f1_suffb}
\end{eqnarray}
That is, since $p(T)=p_0$:
\begin{equation}
\pfrac{F_1}{q}(q=q_0,q_0,T)+\pfrac{F_1}{q_0}(q=q_0,q_0,T)=0
\label{eq:per_f1a}
\end{equation}
 Eq. \ref{eq:per_f1a} defines $n$ equations with $n+1$ variables, $(q_0,T)$, it is an
 under-determined system, and hence we often focus on one of the two following problems:
\begin{enumerate}
\item \titleone: Given a point in the position space $q_0$,
find all periodic orbits going through this point and their associated momentum. Eq.
\ref{eq:per_f1a} defines $n$ equations of a single variable $T$. Taking the norm of the left hand
side yields:
\begin{equation}
\|\pfrac{F_1}{q}(q=q_0,q_0,T)+\pfrac{F_1}{q_0}(q=q_0,q_0,T)\|=0
\label{eq:per_f1b}
\end{equation}
Eq. \ref{eq:per_f1b} is a single equation of one variable that can be solved graphically. To find
the corresponding momentum, we can use either Eq. \ref{eq:f2_in_generala} or Eq.
\ref{eq:f2_in_generalb}:
\begin{eqnarray}
p_0&=&-\pfrac{F_1}{q_0}(q=q_0,q_0,T)\\
p&=&\pfrac{F_1}{q}(q=q_0,q_0,T)
\end{eqnarray}
Both equations provide the same momentum since Eq. \ref{eq:per_f1b} is equivalent to $\|p-p_0\|=0$
and is satisfied.

\item \titlethree: Find all periodic orbits of a given period. Eq. \ref{eq:per_f1a} reduces to a system of $n$
equations with $n$ unknowns, $q_0$. For dynamical systems with $n$ degrees of freedom the solution
may be represented on a $n$-dimensional plot. In practice, solving this problem graphically is
efficient only for systems with at most $3$ degrees of freedom. For Hamiltonian systems with more
than $3$ degrees of freedom, Newton iteration or an equivalent method can be used. When a solution
to Eq.
\ref{eq:per_f1a} is obtained, then we use Eq. \ref{eq:f1_in_generala} or
\ref{eq:f1_in_generalb} to find the corresponding momentum:
\begin{eqnarray}
p_0&=&-\pfrac{F_1}{q_0}(q=q_0,q_0,T)\\
p&=&\pfrac{F_1}{q}(q=q_0,q_0,T)
\end{eqnarray}
\end{enumerate}

\paragraph{Generating function of the second kind}
The search for periodic orbits can also be solved using a generating function of the second kind.
The main difference with the use of $F_1$ is that the system of equations we need to solve does not
reduce to a system of $n$ equations and $n$ functions evaluations (we must solve $2n$ equations).

The generating function $F_2$ allows us to solve a two-point boundary value problem for which the
initial momentum and the position at time $T$ are given. The solution to this problem is found
using Eqns.
\ref{eq:f2_in_generala} and
\ref{eq:f2_in_generalb}.
\begin{eqnarray}
p&=&\pfrac{F_2}{q}(q,p_0,T)\\
q_0&=&\pfrac{F_2}{p_0}(q,p_0,T)
\end{eqnarray}

On the other hand, the boundary value problem is defined by equations \ref{eq:bound_eqa} and
\ref{eq:bound_eqb}. Combining these four equations yields:
\begin{eqnarray}
p_0&=&p(T)\nonumber\\
&=&\pfrac{F_2}{q}(q,p_0,T) \label{eq:per_f2a} \\
q(T)&=&q_0\nonumber\\
&=&\pfrac{F_2}{p_0}(q,p_0,T) \label{eq:per_f2b}
\end{eqnarray}

The system of equations \ref{eq:per_f2a} and \ref{eq:per_f2b} contains $2n$ equations with $2n+1$
variables, and therefore is under-determined. As with $F_1$, we consider two main problems, we
either set the time period or $n$ coordinates of the point in the phase space.

\subsubsection{Examples}
To illustrate the theory developed above, let us consider the Hill three-body problem and let us
find periodic orbits about the Libration point $L_2$ using the generating function of the first
kind $F_1$. To compute $F_1$, we use the algorithm developed by Guibout and Scheeres \cite{gui02}
that computes the Taylor series expansion of the generating functions about a given trajectory,
called the reference trajectory. In this example the reference trajectory is the equilibrium point
$L_2$ and we compute the Taylor series up to order $6$. Since we are working with series expansion,
we will only find periodic orbits that stay within the radius of convergence of the series, not all
periodic orbits.

\begin{figure*}[!h]
\begin{center}
\includegraphics[width=0.7\linewidth]{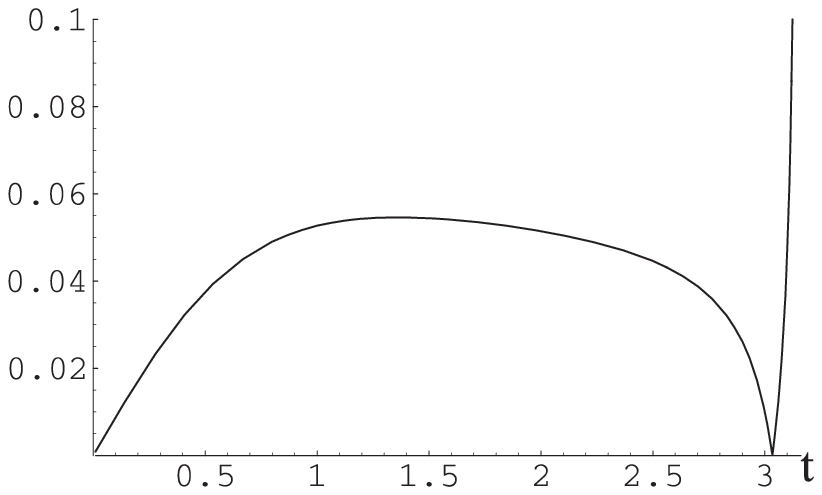}
\caption{\label{fig:f1_norm_nonlin} Plot of $\|\pfrac{F_1}{q}(q=q_0,q_0,T)+\pfrac{F_1}{q_0}(q=q_0,q_0,T)\|$
where $q_0=(0.01,0)$}
\end{center}
\end{figure*}

\paragraph{Search in time domain:} Find all periodic orbits going through the point\footnote{We use
normalized units, for the Sun-Earth-spacecraft system $0.01$ units of length represents about
$21,500km$} $(0.01,0)$. We have seen that this problem can be handled using Eq. \ref{eq:per_f1b}
which is one equation with one variable, $T$. In Figure \ref{fig:f1_norm_nonlin} we have plotted
the left-hand side of Eq. \ref{eq:per_f1b} as a function of time, we obtain a continuous curve
whose points have a particular significance. Let $x$ be a point on that curve whose coordinates are
$x=(t_x,\hp)$. The trajectory whose initial conditions are $q_0=(0.01,0)$,
$p_0=-\pfrac{F_1}{q_0}(q_0,q_0,t_x)$ comes back to its initial position after a time $t_x$ but the
norm of the difference between its initial momentum and its momentum at time $t_x$ is $\hp$. Hence,
any point on the curve whose coordinates are $(t_x,0)$ represents a periodic orbit (not only the
trajectory comes back to its initial position at $t_x$ but the norm of the difference between the
momenta at initial time and at $t_x$ is zero, i.e., the trajectory comes back to its initial state
at $t_x$). In figure \ref{fig:f1_norm_nonlin}, we observe that there exists a periodic orbit of
period $T=3.03353$ going through the point $(0.01,0)$. The corresponding momenta is found using
either Eq. \ref{eq:f2_in_generala} or Eq. \ref{eq:f2_in_generalb} and is $p_0=p=(0,-0.0573157)$.

\begin{figure*}[!h]
\begin{center}
\subfigure[Plot of the solution to the first equation defined by Eq. \ref{eq:per_f1a}]
{\includegraphics[width=0.4\linewidth]{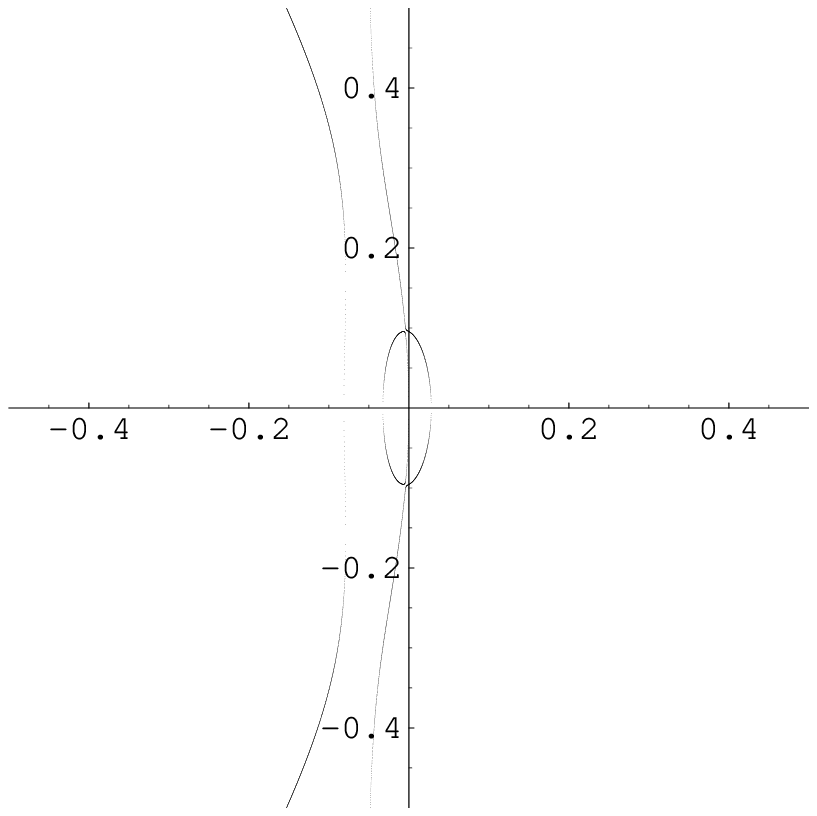}}
\subfigure[Plot of the solution to the second equation defined by Eq. \ref{eq:per_f1a}]
{\includegraphics[width=0.4\linewidth]{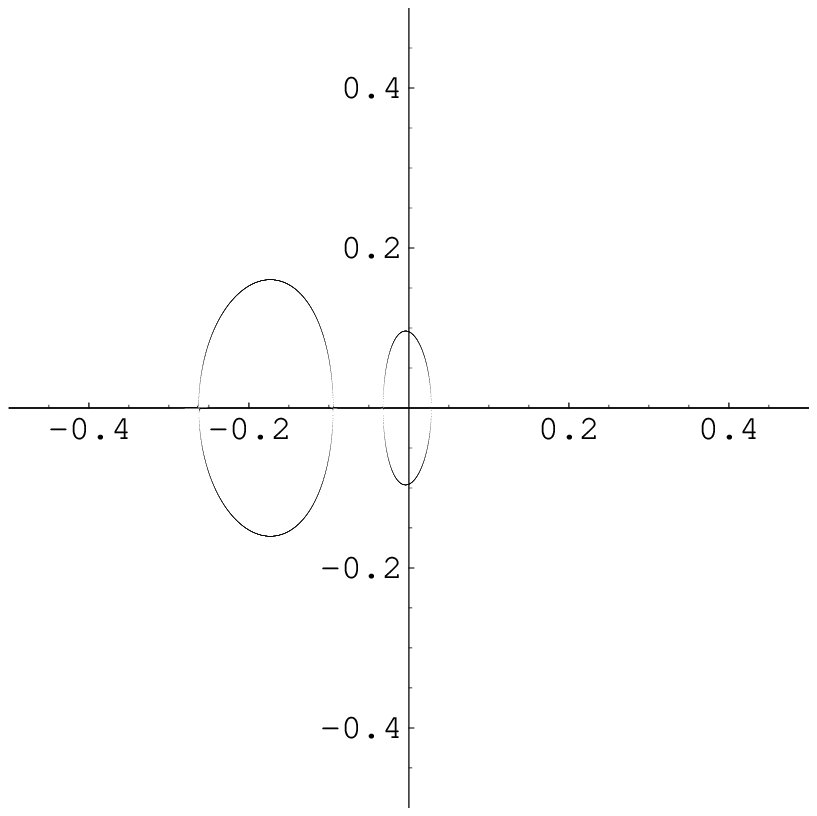}}
\subfigure[Superposition of the two sets of solutions ]
{\includegraphics[width=0.4\linewidth]{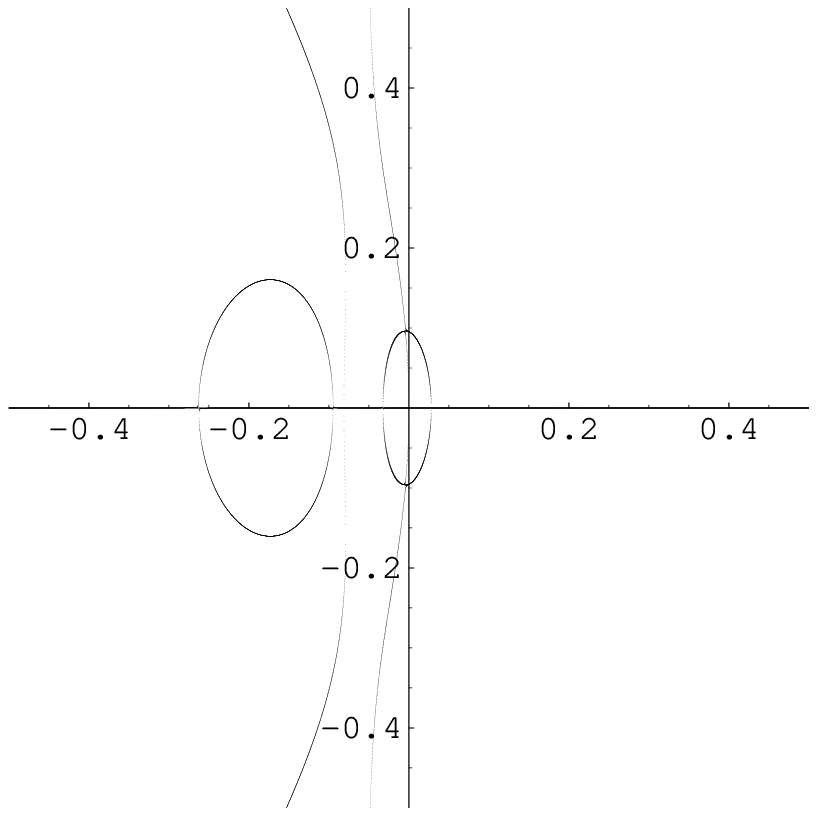}}
\caption{\label{fig:f1_geom_nonlin} Periodic orbits for the nonlinear motion about a Libration point}
\end{center}
\end{figure*}

\begin{figure*}[!h]
\begin{center}
\subfigure[Plot of the solution to the first equation defined by Eq. \ref{eq:per_f1a} for $t=3.033+0.0005n \quad n\in \{1 \cdots 10\}$]
{\includegraphics[width=0.4\linewidth]{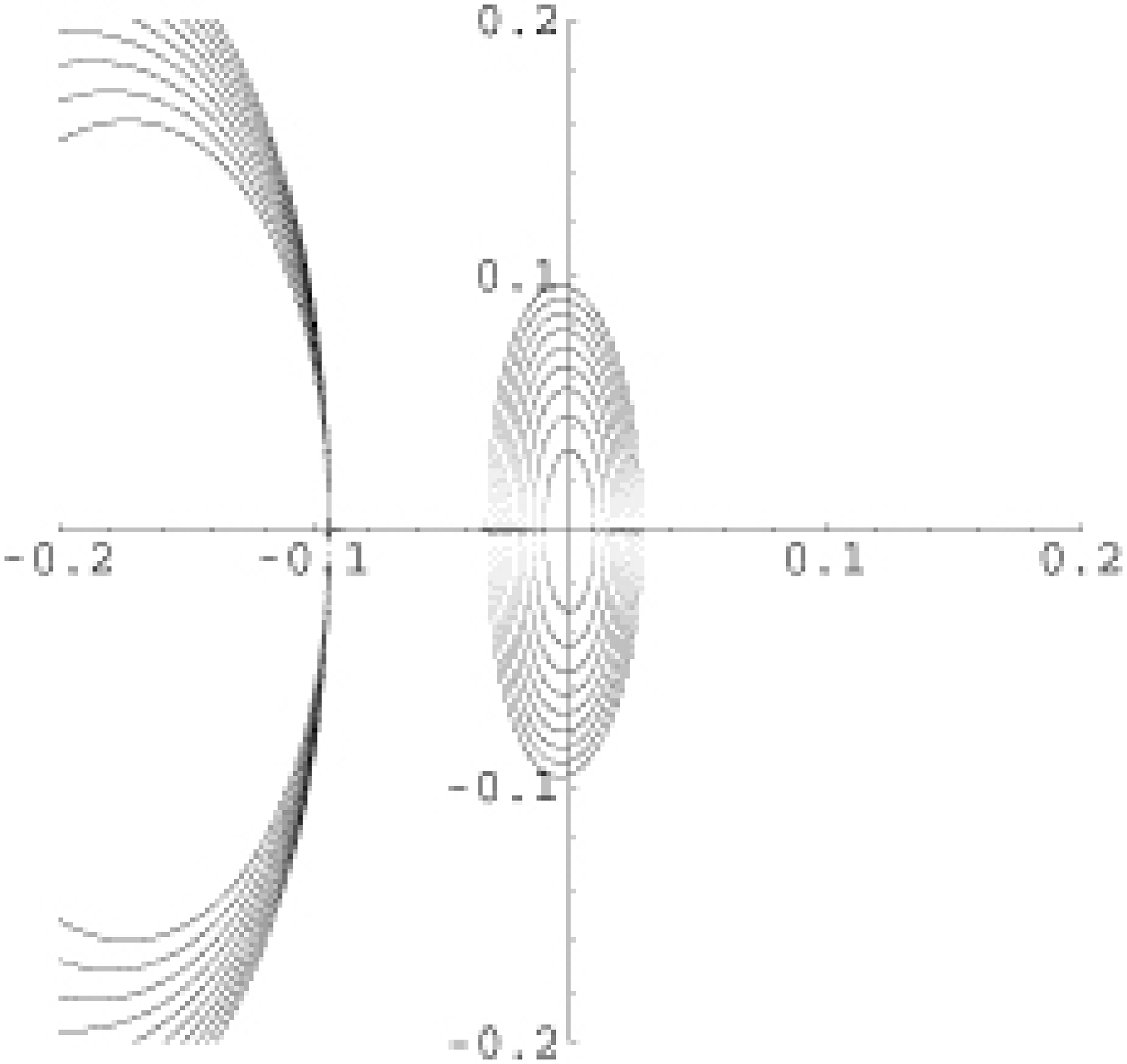}}
\subfigure[Plot of the solution to the second equation defined by Eq. \ref{eq:per_f1a} for $t=3.033+0.0005n \quad n\in \{1 \cdots 10\}$]
{\includegraphics[width=0.4\linewidth]{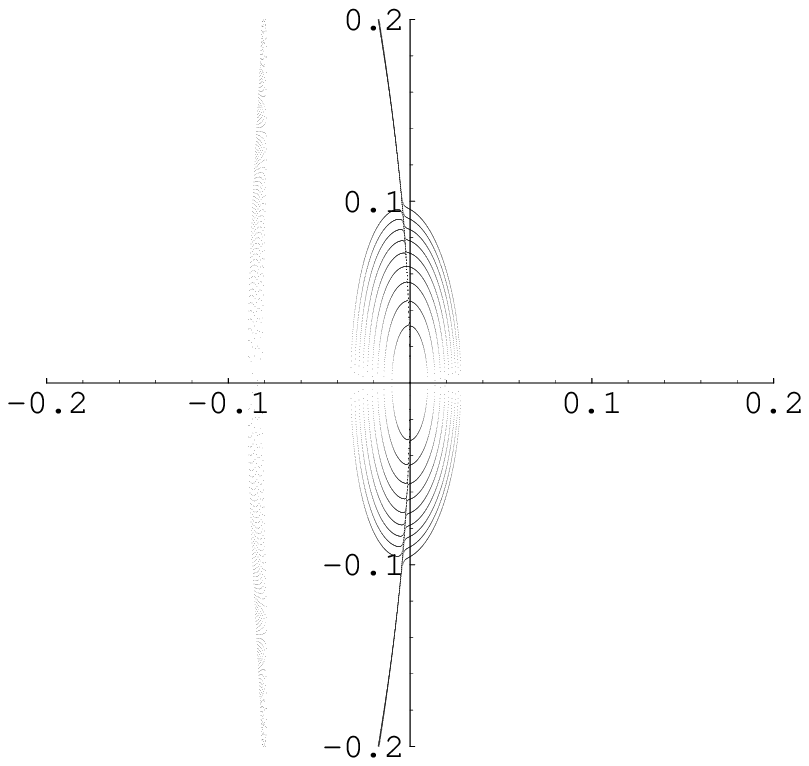}}
\subfigure[Superposition of the two sets of solutions for $t=3.033+0.0005n \quad n\in \{1 \cdots 10\}$]
{\includegraphics[width=0.7\linewidth]{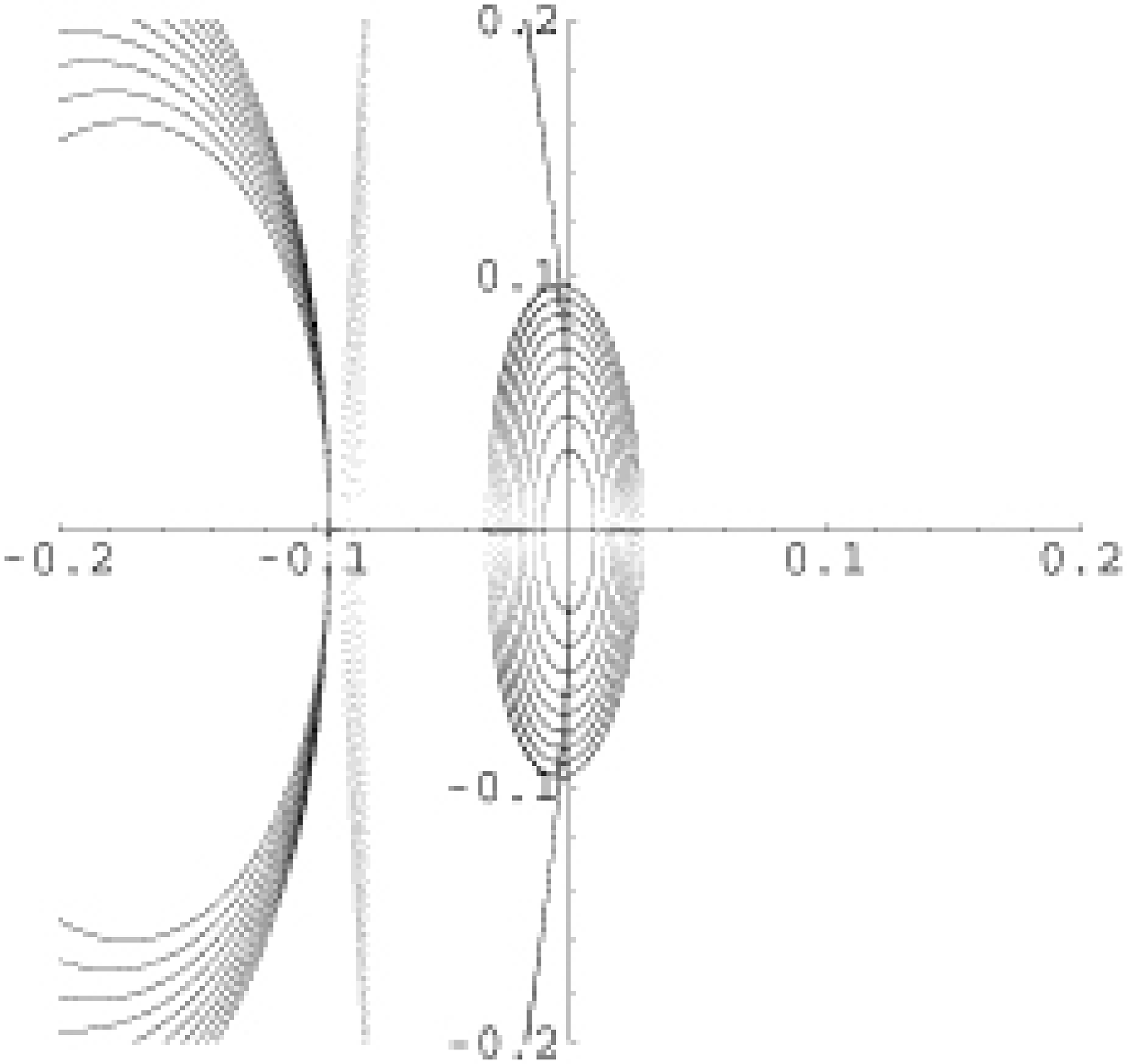}}
\caption{\label{fig:f1_geom_nonlin_multiT} Periodic orbits for the nonlinear motion about a Libration point}
\end{center}
\end{figure*}

\paragraph{Search in position space:} Find all periodic orbits of period $T=3.0345$. To solve this
problem we use Eq. \ref{eq:per_f1a}, which is a system of two equations with two variables
$(q_{0_x},q_{0_y})$. In Fig. \ref{fig:f1_geom_nonlin} we plot solutions to each of these two
equations and then superimpose them to find their intersection, which is the solution to Eq.
\ref{eq:per_f1a}. The solution is a closed curve, i.e., a periodic orbit of the given period.
By plotting the solutions to Eq. \ref{eq:per_f1a} for different periods, we generate a map of a
family of periodic orbits around the Libration point. In Figure
\ref{fig:f1_geom_nonlin_multiT} we plot the solutions to Eq. \ref{eq:per_f1a} for $t=3.033+0.0005n,
 \quad n\in\{0,\cdots,9\}$.

\subsection{Study of equilibrium points}
The generating functions can also be used to study properties of equilibrium points of an
Hamiltonian dynamical system. First, we have proved the equivalence between the state transition
matrix and the generating functions describing relative motion in linear system theory, therefore,
linear terms in the Taylor series expansion of the generating functions about the equilibrium point
provide information on the characteristic time and stability as does the state transition matrix.
The other terms can be used to study the geometry of center, stable and unstable manifolds far from
the equilibrium points where the linear approximation does not hold anymore (but within the radius
of convergence of the Taylor series). The study of center manifolds is a direct application of the
previous section as is readily seen from the example we provided. To find stable and unstable
manifolds we propose a technique that uses generating functions to solve initial value problems,
not two-point boundary value problems. Historically, generating functions were introduced by Jacobi
and used thereafter to solve initial value problems, hence the following technique is not new. We
mention it to show that one is able to fully describe an equilibrium point with only knowledge of
the generating functions.

The idea is to propagate the trajectory of a point that is ``close" to the equilibrium point and on
the linear approximation of the stable (unstable) manifold. Even though this method to find
unstable and stable manifolds is not exact, it is fairly accurate and often used. We then reduce
the search for hyperbolic manifolds to an initial value problem that can be solved using any
generating functions. For simplicity let us consider $F_2$. At the linear level, a point on the
unstable (stable) manifold has coordinates $(q_0,p_0)=(\alpha \hat{u},\alpha \lambda
\hat{u})$ where $\alpha
\ll 1$, $\lambda$ is the characteristic exponent and $\hat{u}$ is the eigenvector defining the
unstable (stable) manifold. Eq. \ref{eq:f2_in_generalb} defines $q(t)$ implicitly:
$$\alpha \hat{u}=q_0=\pfrac{F_2}{p_0}(q,\alpha \lambda \hat{u},t)$$
Once $q(t)$ is found, we find $p(t)$ from Eq. \ref{eq:f2_in_generala}. As $t$ varies, $(q(t),p(t))$
describes the hyperbolic manifolds.

\subsection{Design of spacecraft formation trajectories}
The last application we present concerns the design of a formation of spacecraft. This is again a
direct application of the theory developed in this paper, first introduced by Guibout and
Scheeres\cite{gui02}. This application relies on the fact that the relative dynamics of two
particles evolving in a Hamiltonian dynamical system is Hamiltonian, hence the Hamilton-Jacobi
theory is applicable. To illustrate the use of  generating functions, let us study an example. We
consider a constellation of spacecraft located at the Libration point $L_2$ of Hill's three-body
problem. At a later time $t=t_f$, we want the spacecraft to lie on a circle surrounding the
libration point at a distance of $108,000 km$. What initial velocity is required to transfer to
this circle in time $t_f$, and what will the final velocity be once we arrive at the circle? The
answer will depend, of course, on where we arrive on the circle. In general, this problem must be
solved repeatedly for each point on the circle we wish to transfer to and each transfer time. In
our example we only need to compute the generating functions $F_1$ to be able to compute the answer
as an analytic function of the final location. The method to solve this problem proceeds as
follows: We first compute $F_1$ then we compute the solution to the problem of transferring from
$L_2$ to a point on the final circle where $2$ parameters may vary, the transfer time and the
location on the circle. Then we look at solutions which minimize the total fuel cost of the
maneuver, that is, which minimize the sum of the norm of the initial momentum and the norm of the
final momentum, $\sqrt{|\kpz|^2+|\hp|^2}$. We assume zero momentum in the Hill's rotating frame at
the beginning and end of the maneuver. While not a realistic maneuver, we can use it to exhibit the
applicability of our approach.

\begin{figure*}[!h]
\begin{center}
\includegraphics[width=\linewidth]{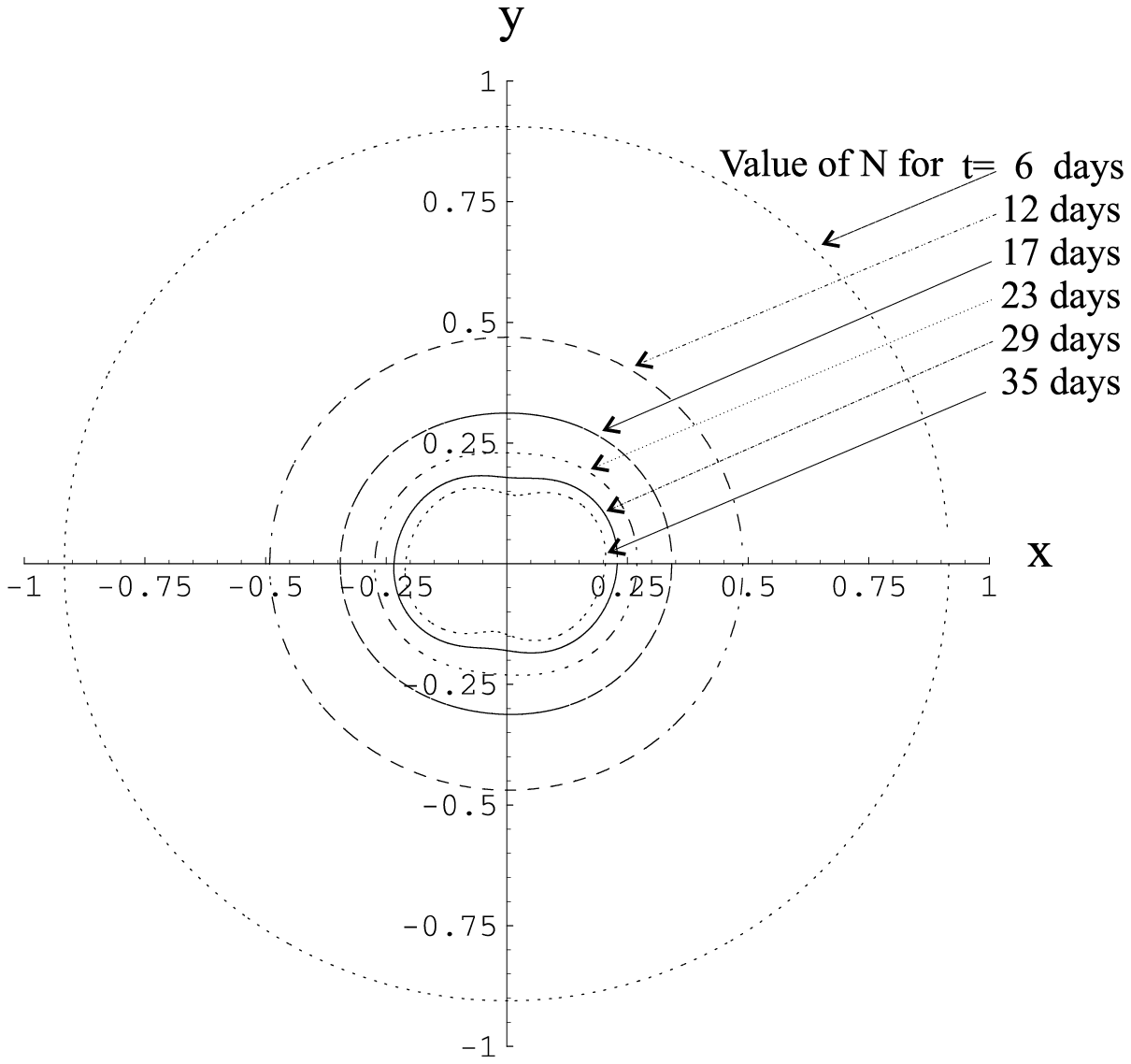}
\caption{\label{fig:optimization1}$\sqrt{|\kpz|^2+|\hp|^2} \hat{q}$ for $t\in [6 days, 35 days ]$}$1
unit \longleftrightarrow 432 m.s^{-1}$
\end{center}
\end{figure*}

\begin{figure*}[!h]
\begin{center}
\includegraphics[width=\linewidth]{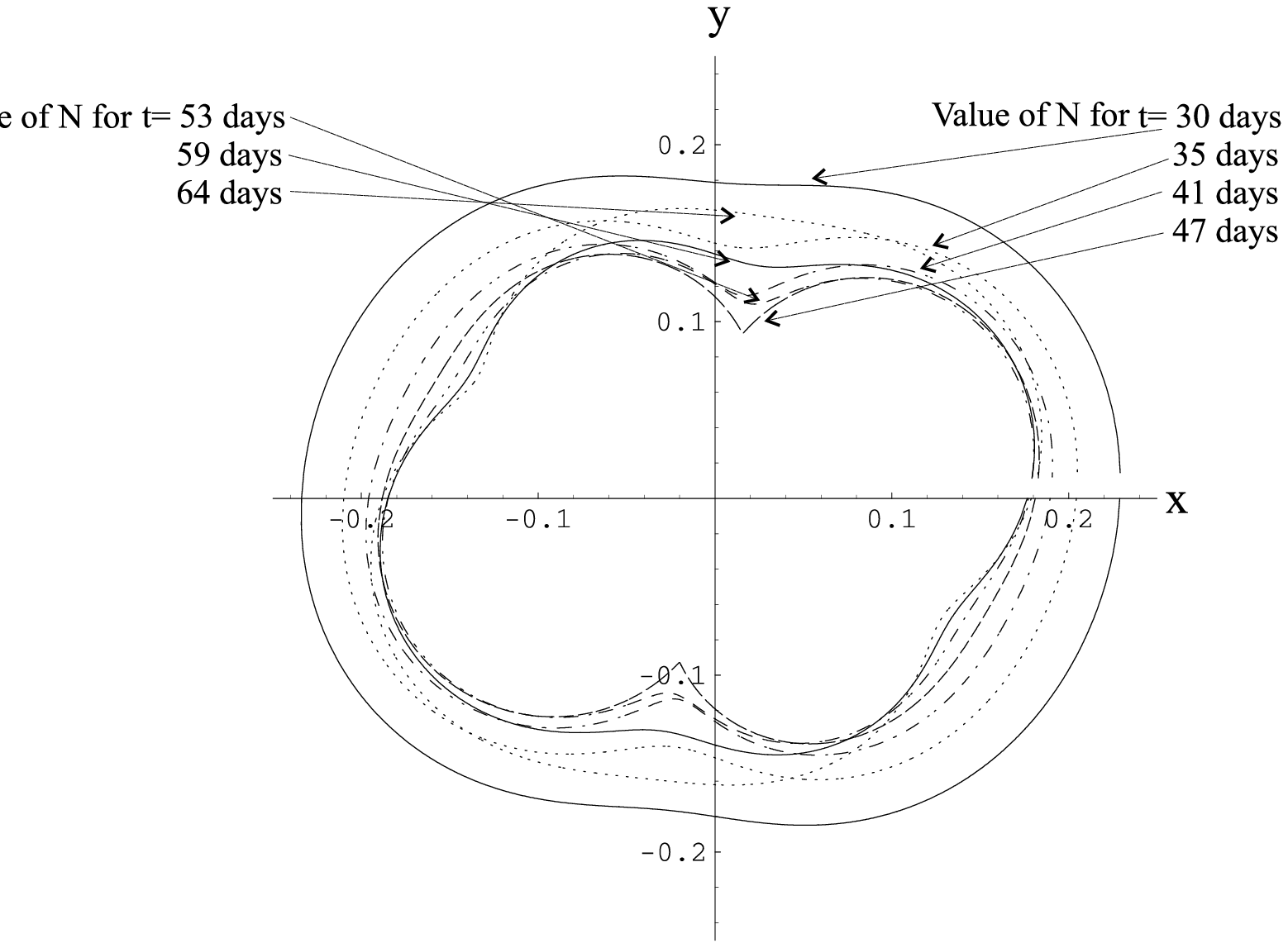}
\caption{\label{fig:optimization2}$\sqrt{|\kpz|^2+|\hp|^2} \hat{q}$ for $t\in [30 days, 64 days ]$}$1
unit \longleftrightarrow 432 m.s^{-1}$
\end{center}
\end{figure*}

\begin{figure*}[!h]
\begin{center}
\includegraphics[width=\linewidth]{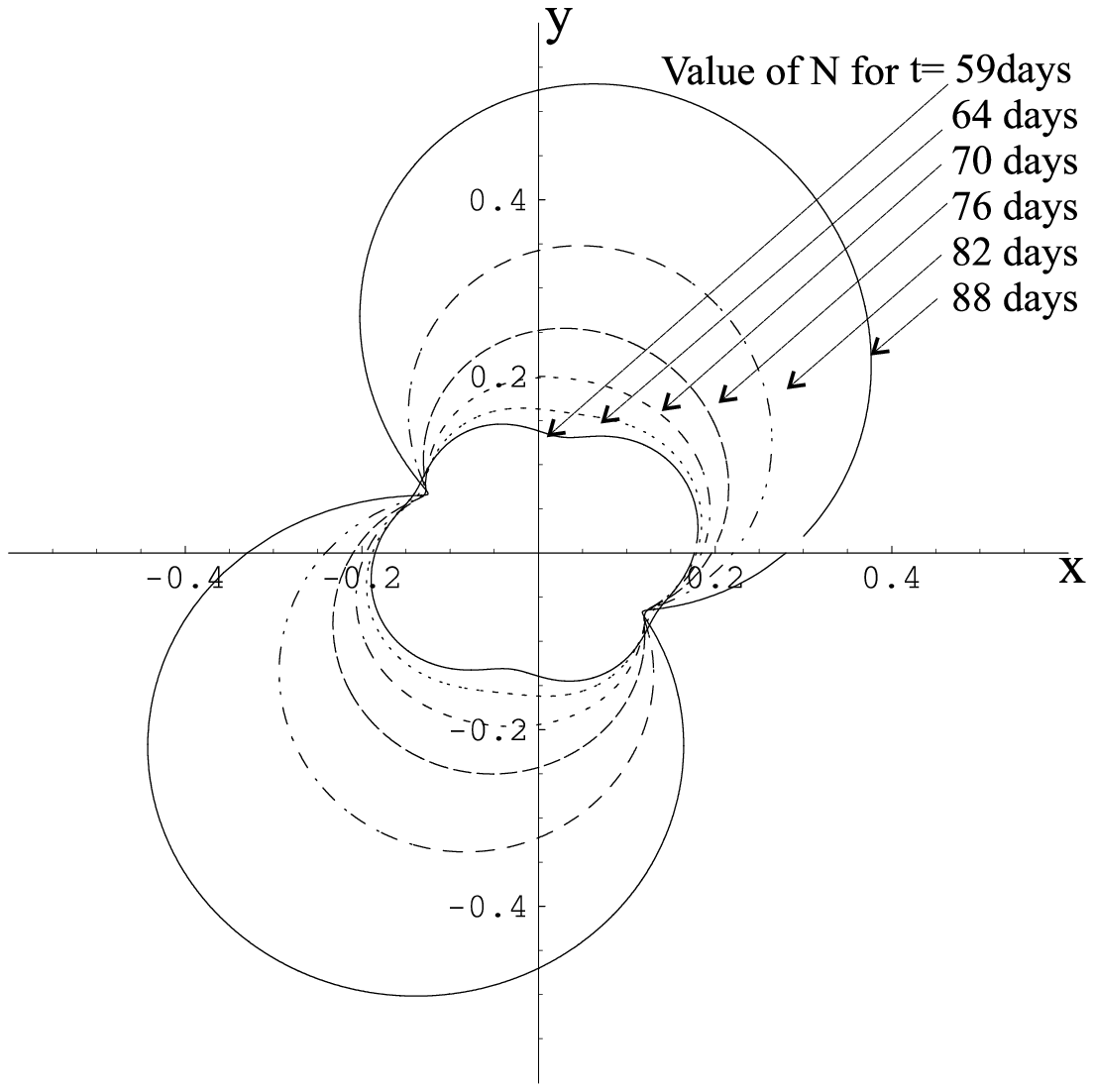}
\caption{\label{fig:optimization3}$\sqrt{|\kpz|^2+|\hp|^2}  \hat{q}$ for $t\in [59 days, 88 days ]$}$1
unit \longleftrightarrow 432 m.s^{-1}$
\end{center}
\end{figure*}

Figures \ref{fig:optimization1}, \ref{fig:optimization2} and \ref{fig:optimization3} show the value
of $\sqrt{|\kpz|^2+|\hp|^2}$ as a function of  position in the final formation at different times
\footnote{Define the final position of the spacecraft as $\hq=\hq \hat{q}$ where $\hq=108,000 km$
and $\hat{q}$ is the unit vector pointing towards the location of the final circle. Then, figures
\ref{fig:optimization1}-\ref{fig:optimization3} represent $\sqrt{|\kpz|^2+|\hp|^2}
\hat{q}$}. We notice three tendencies:
\begin{enumerate}
\item For $t$ less than the characteristic time,
 no matter which  direction the spacecraft leaves $L_2$, it costs essentially the same amount of
fuel to reach the final position and stop (figure \ref{fig:optimization1}).
\item For $t$  larger than the characteristic time, but less than
$47$ days the curve describing
\newline
$\sqrt{|\kpz|^2+|\hp|^2}$ shrinks along a direction $80^\circ$ from the $x$-direction. Thus,
placing a spacecraft on the final circle at an angle of $80^\circ$ or $260^\circ$ from the $x$-axis
provides the lowest cost in fuel (figure \ref{fig:optimization2}).
\item For $t$ larger than $47$ days, the curve describing
$\sqrt{|\kpz|^2+|\hp|^2}$ shrinks along a  direction perpendicular to the previous one, at an angle
of $\sim 170^\circ$ with the x-axis and expands along the $80^\circ$ direction. Thus, there exists
an epoch for  which placing a spacecraft on the final circle at an angle of $170^\circ$ or
$350^\circ$ from the $x$-axis provides the lowest final cost, this happens for $t= 88$ days
(figures \ref{fig:optimization2} and \ref{fig:optimization3}).
\end{enumerate}

To conclude, we see the optimal transfer time to the final circle changes as a function of location
on the circle. While this is to be expected, our results provide direct solutions for this
non-linear boundary value problem.

We now make a few additional remarks to emphasize the advantage of our method. First, additional
spacecraft do not require any additional computations. Hence, our method to design optimal
reconfiguration is valid for infinitely many spacecraft in formation. Second, now that we have
computed the generating functions around the libration point, we are able to analyze any
reconfiguration around the libration point at the cost of evaluating a polynomial
function\footnote{This is especially valuable for missions involving spacecraft that stay close to
$L_2$ since the generating functions in this region can be computed during mission planning. Then
any targeting problem or reconfiguration design can be achieved at the cost of a function
evaluation}. Finally, if the formation of spacecraft is evolving around a base which is on a given
trajectory, we can linearize about this trajectory, and then proceed as in the above examples to
study the reconfiguration problem.

\section*{Conclusions}
This paper describes a novel application of Hamilton-Jacobi theory. We are able to formally solve
any nonlinear two-point boundary value problem using generating functions for the canonical
transformation induced by the phase flow. Many applications of this method are possible, and we
have introduced a few of them, and implemented them successfully. Nevertheless, the method we
propose is based on our ability to obtain generating functions, that is to solve the
Hamilton-Jacobi equation. In general such a solution cannot be found, but for a certain class of
problem an algorithm has been developed \cite{gui02} that converges locally in phase space. A
typical use of this algorithm would be to study the optimal control problem about a known
trajectory, to find families of periodic orbits about an equilibrium point or in the vicinity of
another periodic orbit, and to study spacecraft formation trajectories.

\section*{Appendix I: The circular restricted three-body problem and Hill's three-body problem}
The circular restricted three-body problem is  a three-body problem where the second body is in
circular orbit about the first one and the third body has negligible mass\cite{arn88c}. The
coordinate system is centered at the center of mass of the two bodies with mass and the Hamiltonian
function describing the dynamics of the third body is:

\begin{equation}
H(q_x,q_y,p_x,p_y)=\undemi (p_x^2+p_y^2)+p_x q_y-q_x p_y-
\frac{1-\mu}{\sqrt{(q_x+\mu)^2+q_y^2}}-\frac{\mu}{\sqrt{(q_x-1+\mu)^2+q_y^2}}
\end{equation}
where $q_x=x$, $q_y=y$, $p_x=\dot{x}-y$ and $p_y=\dot{y}+x$. Equations of motion for the third body
can be found from Hamilton's equations:
\begin{eqnarray}
\ddot{x}-2\dot{y}&=&x-(1-\mu)\frac{x+\mu}{((q_x+\mu)^2+q_y^2)^{3/2}}-\mu
\frac{x-1+\mu}{((q_x-1+\mu)^2+q_y^2)^{3/2}}\\
\ddot{y}+2\dot{x}&=&y-(1-\mu)\frac{y}{((q_x+\mu)^2+q_y^2)^{3/2}}-\mu \frac{y}{((q_x-1+\mu)^2+q_y^2)^{3/2}}
\end{eqnarray}
There are five equilibrium points for this system, called the Libration points. $L_2$ is the one
whose coordinates are $( 1.01007,0)$ for a value of $\mu=3.03591\cdot 10^{-6}$.

If the first body has a larger mass than the second one we can expand the equations of motion about
$\mu=0$. Then, shifting the coordinate system so that its center is the second body yields Hill's
formulation of the three-body problem. The equations are motion are:
\begin{eqnarray}
\ddot{x}-2 \dot{y}&=&-\frac{x}{r^3}+3x\\
\ddot{y}+2 \dot{x}&=&-\frac{y}{r^3}\\
\end{eqnarray}
where $r^2=x^2+y^2$.

The Lagrangian then reads:
\begin{equation}
L(q,\dot{q},t)=\undemi (\dot{q}_x^2+\dot{q}_y^2)+\frac{1}{\sqrt{q_x^2+q_y^2}}+\frac{3}{2} q_x^2 -
(\dot{q}_x q_y -\dot{q}_y q_x)
\label{eq:eq1}
\end{equation}

Hence,
\begin{eqnarray}
p_x&=&\pfrac{L}{\dot{q}_x}\nonumber\\
&=&\dot{q}_x-q_y \label{eq:eq2}\\
p_y&=&\pfrac{L}{\dot{q}_y}\nonumber\\
&=&\dot{q}_y+q_x \label{eq:eq3}
\end{eqnarray}

From Eqns. \ref{eq:eq1}, \ref{eq:eq2} and \ref{eq:eq3} we obtain the Hamiltonian function $H$:
\begin{eqnarray}
H(q,p)&=&p_x \dot{q}_x+ p_y \dot{q}_y - L\nonumber\\
&=&\undemi (p_x^2+p_y^2)+(q_y p_x-q_x p_y)-\frac{1}{\sqrt{q_x^2+q_y^2}}+\undemi (q_y^2-2q_x^2)
\end{eqnarray}

There are two equilibrium points for this system, called libration points.  Their coordinates are
$L_1(-\left(\frac{1}{3}\right)^{1/3},0)$ and $L_2(\left(\frac{1}{3}\right)^{1/3},0)$

\newpage

\newpage

\newpage

\newpage

\newpage

\newpage

\newpage

\bibliographystyle{plain}
\bibliography{./submitted}

\end{document}